\documentclass[11pt]{article}
\usepackage{mathrsfs}
\usepackage{amssymb}
\usepackage{amsmath}
\usepackage{amsbsy}
\usepackage{epsfig}
\usepackage{enumerate}
\usepackage{bm}
\usepackage{color}

\newcommand{\blue}[1]{\textcolor{blue}{#1}}

\usepackage[colorlinks,linkcolor=blue]{hyperref}

\newcommand{\ARXIV}[1]{\href{http://arXiv.org/abs/#1}{\blue{arXiv:#1}}}

\newtheorem{lemma}{Lemma}[section]
\newtheorem{theorem}{Theorem}[section]
\newtheorem{proposition}{Proposition}[section]

\newtheorem{corollary}{Corollary}[section]
\newtheorem{definition}{Definition}[section]
\newtheorem{remark}{Remark}[section]

\def\pr{\textsf{P}}

\def\Sbep{\widehat{\mathbb E}}
\def\cSbep{\widehat{\mathcal E}}
\def\vSbep{\breve{\mathbb E}}
\def\vcSbep{\breve{\mathcal E}}
\def\Capc{\mathbb V}
\def\cCapc{\mathcal V}
\def\upCapc{\widehat{\mathbb V}}
\def\lowCapc{\widehat{\mathcal V}}
\def\outCapc{\widehat{\mathbb V}^{\ast}}
\def\outcCapc{\widehat{\mathcal V}^{\ast}}

\renewcommand{\baselinestretch}{1.5}

\begin{document}

\begin{center}
{\LARGE\bf The Limit Points of the Strong Law of Large Numbers under Sub-linear Expectations}\\
\ARXIV{2311.11100}
\end{center}

\begin{center}
{\sc Li-Xin ZHANG\footnote{This work was supported by grants from the National Key R\&D Program of China (No. 2024YFA1013502), the NSF of China (Grant No. U23A2064), and the Summit Advancement Disciplines of Zhejiang Province (Zhejiang Gongshang University--Statistics).}}
\\
{\sl \small School of Statistics and Data Science, Zhejiang Gongshang University, Hangzhou 310018}\\
(Email: stazlx@mail.zjgsu.edu.cn)\\
\end{center}

\renewcommand{\abstractname}{~}
\begin{abstract}
{\bf Abstract:} Let $\{X_n;n\ge 1\}$ be a sequence of independent and identically distributed random variables relative to a regular sub-linear expectation $(\Omega,\mathscr{H},\widehat{\mathbb E})$  with  upper mean $\overline{\mu}$,  lower mean $\underline{\mu}$,  and finite Choquet expectation or $\lim_{c\to\infty} \widehat{\mathbb E}[(|X_1|-c)^+]=0$. Then for any Borel-measurable function $\varphi(x_1,\ldots,x_d)$ on $\mathbb R^d$ or continuous function $\varphi(x_1,x_2,\ldots)$ on $\mathbb R^{\mathbb N}$,
$\sum_{i=1}^nX_i/n$ converges to $\underline{\mu}\wedge \varphi(X_1,X_2,\ldots)\wedge \overline{\mu}$ with upper capacity $1$.
The set of limit points of $\sum_{i=1}^nX_i/n$ is the whole mean interval $[\underline{\mu},\overline{\mu}]$ with upper capacity $1$, as is any random subset of the mean interval whose boundaries are continuous functions or finite-dimensional Borel-measurable functions of $(X_1, X_2, \ldots)$, also with upper capacity $1$.



\vspace{-3mm}
\end{abstract}

\renewcommand{\baselinestretch}{1.2}

\section{Introduction and Notations}\label{sect1}
\setcounter{equation}{0}

Let $\{X_n;n\ge 1\}$ be a sequence of independent and identically distributed random variables (i.i.d.) on a probability space $(\Omega,\mathcal F, \pr)$. Denote $S_n=\sum_{i=1}^n X_i$. Kolmogorov's \cite{Kol} strong law of large numbers (cf. Theorem 3.2.2 of Stout \cite{Stout74}, Theorem 6.11 of Petrov \cite{Petrov95}) states that
\begin{equation}\label{eqKol1}
\pr\left(\lim_{n\to \infty}\frac{S_n}{n}=b\right)=1
\end{equation}
if and only if
\begin{equation}\label{eqKol2}
E_{\pr}[|X_1|]<\infty \quad \text{and}\quad E_{\pr}[X_1]=b,
\end{equation}
where $E_{\pr}$ is the expectation with respect to the probability measure $\pr$. When the probability measure $\pr$ is uncertain, one may consider a family $\mathscr{P}$ of probability measures and define $\Sbep[X]=\sup_{P\in \mathscr{P}}E_P[X]$. Then $\Sbep$ is no longer a linear expectation, and the related capacities $\cCapc(A)=\inf_{P\in\mathscr{P}}P(A)$ and $\Capc(A)=\sup_{P\in \mathscr{P}}P(A)$ are not additive. Instead, $\Sbep$ is sub-linear in the sense that $\Sbep[aX+bY]\le a\Sbep[X]+b\Sbep[Y]$ if $a,b\ge 0$, and $\Capc(A)=\sup_{P\in \mathscr{P}}P(A)$ is sub-additive. Limit laws for non-additive capacities and non-linear expectations have attracted increasing interest. Marinacci \cite{Mar1999}, Maccheroni and Marinacci \cite{MacMar2005}, and Ter\'an \cite{Teran14} obtained the strong law of large numbers for a sequence of random variables $\{X_n;n\ge 1\}$ which are independent (preindependent, pairwise preindependent) and identically distributed relative to a {\em completely monotone} and {\em continuous} (or topologically continuous) capacity $\cCapc$ by showing
\begin{equation}
\label{eq:LLNCapc}
\cCapc\left(\int X_1d\cCapc \le \liminf_{n\to\infty}\frac{S_n}{n}\le \limsup_{n\to\infty}\frac{S_n}{n}\le -\int(-X_1)d\cCapc\right)=1.
\end{equation}
The limit values are points in an interval with separate bounds instead of a single limit. In \eqref{eq:LLNCapc}, $\int Xd\nu$ is defined to be the Choquet integral of $X$ with respect to $\nu$, and $-\int(-X)d\nu=\int Xd\overline{\nu}$, $\overline{\nu}(A)=1-\nu(A^c)$. A natural question is whether the quality of the bounds in \eqref{eq:LLNCapc} holds or not.

Ter\'an \cite{Teran14} showed by examples that the set of cluster points of $S_n/n$ can be either the whole mean interval $[\int X_1d\cCapc, -\int (-X_1) d\cCapc]$ or a proper subset, and left the following problem as one of his open questions:
\begin{description}
\item[Q1] {\em Does choosing an independence definition using some specific product capacity impose interesting restrictions on the behavior of $S_n/n$?}
\end{description}

For a completely monotone capacity $\nu$, $\int Xd\overline{\nu}$ has the representation
$$ \int Xd\overline{\nu}=\sup_{P\in \mathrm{core}(\nu)}E_P[X], $$
and is a sub-linear expectation, where $\mathrm{core}(\nu)=\{P:P \text{ is a probability measure with } P\ge \nu\}$, cf. Ter\'an \cite[Proposition 3.5]{Teran14}.

Peng \cite{peng2008,peng2019} introduced the concepts of independence and identical distribution expressed through the sub-linear expectation $\Sbep$ instead of a capacity $\cCapc$, and established the weak law of large numbers and central limit theorem for independent and identically distributed random variables. Also, the limit in the weak law of large numbers is described by an interval. Because $\cCapc$ is the infimum of a family of probabilities, the sequence $S_n/n$ cannot be expected to satisfy the simultaneous LLN with a unique limit so that \eqref{eqKol1} holds with $\cCapc$ taking the place of $P$ (cf. Ter\'an \cite{Teran14}). The purpose of this paper is to find conditions and the limit value $b$ for \eqref{eqKol1} to hold for some $P\in \mathscr{P}$. It will be shown that, in contrast to the single constant limit in the classical Kolmogorov strong law of large numbers, for a sequence $\{X_n;n\ge 1\}$ of random variables which are independent and identically distributed relative to a regular sub-linear expectation, any continuous or finite-dimensional Borel-measurable function $f(X_1, X_2,\ldots)$ is a limit value of $\lim_{n\to\infty}\frac{S_n}{n}$ almost surely under some $P\in\mathscr{P}$, as well as of $\liminf_{n\to\infty}\frac{S_n}{n}$ and $\limsup_{n\to\infty}\frac{S_n}{n}$. We will find that if independence is defined relative to the sub-linear expectation, then the answer to Ter\'an's \cite{Teran14} question (Q1) is positive, the mean interval is an optimal limit set, and all random sub-intervals with continuous or finite-dimensional Borel-measurable bounds are also sets of cluster points of $S_n/n$ with full upper capacity. A claim of Peng \cite{peng2009,peng2010}, which was shown to be wrong by Ter\'an \cite{Teran18} via a special counterexample, is proved correct when the sub-linear expectation is regular. The main results will be given in Section \ref{sectMain}. We recall the framework of sub-linear expectations in this section.

We use the framework and notations of Peng \cite{peng2008,peng2019}. Let $(\Omega, \mathcal F)$ be a given measurable space and let $\mathscr{H}$ be a linear space of real functions defined on $(\Omega,\mathcal F)$ such that if $X_1,\ldots, X_n \in \mathscr{H}$, then $\varphi(X_1,\ldots,X_n)\in \mathscr{H}$ for each $\varphi\in C_{l,Lip}(\mathbb R^n)$, where $C_{l,Lip}(\mathbb R^n)$ denotes the linear space of (local Lipschitz) functions $\varphi$ satisfying
\begin{eqnarray*}
& |\varphi(\bm x) - \varphi(\bm y)| \le  C(1 + |\bm x|^m + |\bm y|^m)|\bm x- \bm y|, \;\; \forall \bm x, \bm y \in \mathbb R^n,&\\
& \text{for some }  C > 0, m \in \mathbb  N \text{ depending on } \varphi. &
\end{eqnarray*}
We also denote by $C_{Lip}(\mathbb R^n)$, $C_{b,Lip}(\mathbb R^n)$, $C_b(\mathbb R^n)$, and $C_{b}^{\infty}(\mathbb R^n)$ the spaces of Lipschitz functions, bounded Lipschitz functions, bounded continuous functions, and bounded continuous functions having bounded derivatives of each order on $n$-dimensional real space $\mathbb R^n$, respectively.

\begin{definition}\label{def1.1}
A sub-linear expectation $\Sbep$ on $\mathscr{H}$ is a function $\Sbep: \mathscr{H}\to \overline{\mathbb R}$ satisfying the following properties: for all $X, Y \in \mathscr H$,
\begin{description}
\item[\rm (a)] If $X \ge Y$, then $\Sbep [X]\ge \Sbep [Y]$;
\item[\rm (b)] $\Sbep [c] = c$;
\item[\rm (c)] $\Sbep[X+Y]\le \Sbep [X] +\Sbep [Y ]$ whenever $\Sbep [X] +\Sbep [Y ]$ is not of the form $\infty-\infty$ or $-\infty+\infty$;
\item[\rm (d)] $\Sbep [\lambda X] = \lambda \Sbep  [X]$, $\lambda\ge 0$.
\end{description}
Here $\overline{\mathbb R}=[-\infty, \infty]$. The triple $(\Omega, \mathscr{H}, \Sbep)$ is called a sub-linear expectation space. Given a sub-linear expectation $\Sbep$, the conjugate expectation $\cSbep$ of $\Sbep$ is defined by
$$ \cSbep[X]:=-\Sbep[-X], \;\; \forall X\in \mathscr{H}. $$
\end{definition}

By Theorem 1.2.1 of Peng \cite{peng2019}, there exists a family of finite-additive linear expectations $E_{\theta}: \mathscr{H}\to \overline{\mathbb R}$ indexed by $\theta\in \Theta$, such that
\begin{equation}\label{linearexpression}
\Sbep[X]=\max_{\theta\in \Theta} E_{\theta}[X] \quad \text{for } X \in \mathscr{H}.
\end{equation}
Moreover, for each $X\in \mathscr{H}$, there exists $\theta_X\in \Theta$ such that $\Sbep[X]=E_{\theta_X}[X]$ if $\Sbep[X]$ is finite.

Ter\'an \cite{Teran18} showed that a sub-linear expectation coincides with Walley's \cite{Walley1991} coherent upper prevision. Peng \cite{peng2008,peng2019} introduced concepts of independence, identical distribution, etc., under the sub-linear expectation.

\begin{definition}[{\sl See Peng \cite{peng2008,peng2019}}]
\begin{description}
\item[\rm (i)] ({\em Identical distribution}) Let $\bm X_1$ and $\bm X_2$ be two $n$-dimensional random vectors defined respectively in sub-linear expectation spaces $(\Omega_1, \mathscr{H}_1, \Sbep_1)$ and $(\Omega_2, \mathscr{H}_2, \Sbep_2)$. They are called identically distributed, denoted by $\bm X_1\overset{d}= \bm X_2$, if
$$ \Sbep_1[\varphi(\bm X_1)]=\Sbep_2[\varphi(\bm X_2)], \;\; \forall \varphi\in C_{b,Lip}(\mathbb R^n). $$
A sequence $\{X_n;n\ge 1\}$ of random variables is said to be identically distributed if $X_i\overset{d}= X_1$ for each $i\ge 1$.
\item[\rm (ii)] ({\em Independence}) In a sub-linear expectation space $(\Omega, \mathscr{H}, \Sbep)$, a random vector $\bm Y = (Y_1, \ldots, Y_n)$, $Y_i \in \mathscr{H}$, is said to be independent of another random vector $\bm X = (X_1, \ldots, X_m)$, $X_i \in \mathscr{H}$, under $\Sbep$ if for each test function $\varphi\in C_{b,Lip}(\mathbb R^m \times \mathbb R^n)$ we have
$ \Sbep [\varphi(\bm X, \bm Y )] = \Sbep \big[\Sbep[\varphi(\bm x, \bm Y )]\big|_{\bm x=\bm X}\big]$.

A sequence of random variables $\{X_n; n\ge 1\}$ is said to be independent if $X_{i+1}$ is independent of $(X_{1},\ldots, X_i)$ for each $i\ge 1$.
\end{description}
\end{definition}

Let $(\Omega, \mathscr{H}, \Sbep)$ be a sub-linear expectation space. We denote by $(\upCapc,\lowCapc)$ a pair of capacities defined by
\begin{equation}\label{equpCapc}
\upCapc(A):=\inf\{\Sbep[\xi]: I_A\le \xi, \xi\in\mathscr{H}\},\quad \lowCapc(A)=1-\upCapc(A^c), \quad \forall A\in \mathcal F,
\end{equation}
where $A^c$ is the complement of $A$. Then $\upCapc$ is a sub-additive capacity with the property that
\begin{equation}\label{eq1.4}
\begin{matrix}
\Sbep[f]\le \upCapc(A)\le \Sbep[g]\\
\cSbep[f]\le \lowCapc(A)\le \cSbep[g]
\end{matrix}
\quad
\text{if } 0\le f\le I_A\le g, f,g \in \mathscr{H} \text{ and } A\in \mathcal F.
\end{equation}
We call $\upCapc$ and $\lowCapc$ the upper capacity and the lower capacity, respectively.

Also, we define the Choquet integrals/expectations $(C_{\upCapc},C_{\lowCapc})$ by
$$ C_V[X]=\int_0^{\infty} V(X\ge t)dt +\int_{-\infty}^0\left[V(X\ge t)-1\right]dt $$
with $V$ being replaced by $\upCapc$ and $\lowCapc$, respectively.

If $\Capc$ on the sub-linear expectation space $(\Omega,\mathscr{H},\Sbep)$ and $\widetilde{\Capc}$ on the sub-linear expectation space $(\widetilde{\Omega},\widetilde{\mathscr{H}},\widetilde{\mathbb E})$ are two capacities having property \eqref{eq1.4}, then for any random variables $X\in \mathscr{H}$ and $\tilde X\in \widetilde{\mathscr{H}}$ with $X\overset{d}=\tilde X$, we have
\begin{equation}\label{eqV-V}
\Capc(X\ge x+\epsilon)\le \widetilde{\Capc}(\tilde X\ge x)\le \Capc(X\ge x-\epsilon)\;\; \text{for all } \epsilon>0 \text{ and } x,
\end{equation}
and
\begin{equation}\label{eqV-V3}
C_{\Capc}(X)=C_{\widetilde{\Capc}}(\tilde X).
\end{equation}

Because a capacity $\upCapc$ may not be countably sub-additive, so that the Borel--Cantelli lemma is not valid, we consider its countably sub-additive extension $\outCapc$ defined by
\begin{equation}\label{outcapc}
\outCapc(A)=\inf\Big\{\sum_{n=1}^{\infty}\upCapc(A_n): A\subset \bigcup_{n=1}^{\infty}A_n\Big\},\quad \outcCapc(A)=1-\outCapc(A^c),\quad A\in\mathcal F.
\end{equation}
As shown in Zhang \cite{Zh16a}, $\upCapc$ (resp. $\outCapc$) is the largest sub-additive (resp. countably sub-additive) capacity in the sense that if $V$ is also a sub-additive (resp. countably sub-additive) capacity satisfying $V(A)\le \Sbep[g]$ whenever $I_A\le g\in \mathscr{H}$, then $V(A)\le \Capc(A)$ (resp. $V(A)\le \outCapc(A)$). Also, we denote $\outcCapc(A)=1-\outCapc(A^c)$, still referring to $\outCapc$ and $\outcCapc$ as the upper capacity and lower capacity, respectively.

For real numbers $x$ and $y$, denote $x\vee y=\max(x,y)$, $x\wedge y=\min(x,y)$, $x^+=\max(0,x)$, and $x^-=\max(0,-x)$. For a function $\varphi$, denote $\|\varphi\|_{\infty}= \sup_x|\varphi(x)|$. For a random variable $X$, because $XI\{|X|\le c\}$ may not be in $\mathscr{H}$, we will truncate it in the form $(-c)\vee X\wedge c$, denoted by $X^{(c)}$. Define
$\vSbep[X]=\lim\limits_{c\to\infty}\Sbep[X^{(c)}]$ if the limit exists, and $\vcSbep[X]=-\vSbep[-X]$. It is obvious that $\vSbep[X]$ exists when $C_{\upCapc}(|X_1|)<\infty$ and $|\vSbep[X]|\le \vSbep[|X|]\le C_{\upCapc}(|X_1|)$,
$C_{\lowCapc}(X)\le \vcSbep[X]\le \vSbep[X]\le C_{\upCapc}(X)$, and $\vSbep[(|X|-c)^+]\to 0$ as $c\to\infty$. Further, if $\Sbep[(|X|-c)^+]\to 0$ as $c\to\infty$, then $\vSbep[X]=\Sbep[X]$.

Finally, for a probability measure $P$, we denote $E_P[X]=\int X dP$. For a family of probability measures $\mathscr{P}$, we denote
$$ \Capc^{\mathscr{P}}(A)=\sup_{P\in\mathscr{P}}P(A)\quad \text{and}\quad \mathbb E^{\mathscr{P}}[X]=\sup_{P\in \mathscr{P}}E_P[X]. $$

\section{Main Results}\label{sectMain}
\setcounter{equation}{0}

Let $\{X_n;n\ge 1\}$ be a sequence of i.i.d. random variables in a sub-linear expectation space $(\Omega,\mathscr{H},\Sbep)$.
Denote $\bm X=(X_1,X_2,\ldots)$, $\underline{\mu}=\vcSbep[X]$, $\overline{\mu}=\vSbep[X]$, $S_n=\sum_{i=1}^n X_i$, and $\mathscr{H}_b=\{f\in \mathscr{H}; f \text{ is bounded}\}$. Peng proved the following weak law of large numbers (WLLN) for sub-linear expectations.

\begin{proposition}\label{prop1} (Peng \cite[Theorem 2.4.1]{peng2019})
Suppose
\begin{align}
\label{eq:Pengcondition}
\Sbep[(|X_1|-c)^+]\to 0 \text{ as } c\to \infty.
\end{align}
Then, for any $\varphi\in C_{l,Lip}(\mathbb R)$,
\begin{equation}\label{eq:WLLN}
\Sbep\left[\varphi\left(\frac{S_n}{n}\right)\right]\to \sup_{x\in [\underline{\mu},\overline{\mu}]}\varphi(x).
\end{equation}
\end{proposition}

Based on this result, in \cite[Sect. 5.2]{peng2009} and \cite[Sect. 2.5]{peng2010}, Peng pointed out that, when $n\to \infty$, the number $S_n/n$ can take any value inside $[\underline{\mu},\overline{\mu}]$, and proposed a Monte Carlo method relying on the following formula to obtain an almost sure sample approximation to the expectation:
\begin{equation}\label{eq:PengMC}
\hat{\mathbb M}[\varphi]=:
\limsup_{n\to \infty} \frac{1}{n}\sum_{k=1}^n \psi(X_k)=\Sbep[\psi(X_1)], \;\; \psi\in C_{l,Lip}(\mathbb R).
\end{equation}
However, Ter\'an \cite{Teran18} showed that \eqref{eq:PengMC} and \eqref{eq:sec2.1} below do not hold in general by a counterexample in which $\underline{\mu}=\cSbep[X_1]=0$, $\overline{\mu}=\Sbep[X_1]=1$, and $S_n(\omega)/n\to 0$ for all $\omega\in \Omega$. Based on his example, Ter\'an concluded that ``Peng's generalization of the Monte Carlo method is wrong''. A natural question is that under what kinds of conditions, \eqref{eq:PengMC} remains true. 

Chen \cite{Chen16} and Zhang \cite{Zh16a} obtained the following strong law of large numbers (SLLN) similar to \eqref{eq:LLNCapc}:
\begin{equation}\label{eq:SLLN0}
\outcCapc\left(  \underline{\mu}\le  \liminf_{n\to \infty} \frac{S_n}{n} \le \limsup_{n\to \infty}\frac{S_n}{n} \le  \overline{\mu}\right)=1
\end{equation}
under the conditions $\Sbep[|X_1|^{1+\alpha}]<\infty$ and
\begin{equation}\label{eq:moment}
C_{\upCapc}(|X_1|)<\infty,
\end{equation}
respectively.

To obtain the precise limits so that \eqref{eq:PengMC} holds with $\upCapc$-capacity $1$, or equivalently, 
\begin{equation}\label{eq:sec2.1}
\limsup_{n\to \infty} \frac{S_n}{n}=\overline{\mu},
\end{equation}
\begin{equation}\label{eq:sec2.2}
\liminf_{n\to \infty}\frac{S_n}{n}=\underline{\mu}
\end{equation}
hold with $\upCapc$-capacity $1$, Chen \cite{Chen16} and Zhang \cite{Zh16a} assumed that the sub-additive capacity $\upCapc$ related to the sub-linear expectation $\Sbep$ is continuous from above in the sense that $\upCapc(A_n)\searrow \upCapc(A)$ for all $A_n,A\in \mathcal{F}$ with $A_n\searrow A$.
However, Zhang \cite{ZhangLIL21} pointed out that, under the framework of sub-linear expectations, the continuity of the capacity is a very strict condition by showing that if a capacity $\Capc$ satisfying \eqref{eq1.4} is continuous and $\{Y_i;i\ge 1\}$ is a sequence of i.i.d. random variables in the sub-linear expectation space, then the sub-linear expectation $\Sbep$ is a linear expectation on the space of local Lipschitz functions of $\{Y_i; i\ge 1\}$ and $\Capc$ is a probability measure on $\sigma(Y_i;i\ge 1)$; cf. the proofs of Proposition 4.1, Lemma 4.2 (b), and Lemma 4.3 of Zhang \cite{ZhangLIL21}. To consider sufficient and necessary conditions of the SLLN and the law of the iterated logarithm for i.i.d. random variables, Zhang \cite{ZhangLLN, ZhangLIL21} proposed a reasonable condition on the sub-linear expectation as follows to replace the continuity of the capacity.

\begin{description}
\item[\rm (CC)] The sub-linear expectation $\Sbep$ on $\mathscr{H}_b$ has the representation
\begin{equation}
\label{eq:representbyP}
\Sbep[X]=\sup_{P\in \mathscr{P}}E_P[X], \; X\in \mathscr{H}_b,
\end{equation}
where $\mathscr{P}$ is a countably-dimensionally weakly compact family of probability measures on $\mathscr{H}_b$ in the sense that, for any $Y_1,Y_2,\ldots \in \mathscr{H}_b$, the family $\mathscr{P}\bm Y^{-1}=:\{\overline{P}:
\overline{P}(A)=P(\bm Y\in A), A\in \mathbb R^{\mathbb N}, P\in \mathscr{P}\}$ is a weakly compact family of probability measures on the metric space $\mathbb R^{\mathbb N}=\{\bm x=(x_1,x_2,\ldots): x_i\in\mathbb R, i=1,2,\ldots\}$ with metric $\mathrm{dist}(\bm x,\bm y)=\sum_{i=1}^{\infty}(|x_i-y_i|\wedge 1)/2^i$, where $\bm Y=(Y_1, Y_2,\ldots)$.
\end{description}

In \eqref{eq:representbyP}, $\mathscr{P}$ is called a representation of $\Sbep$. It has been shown that condition (CC) is widely satisfied and is equivalent to $\Sbep$ being regular on $\mathscr{H}_b$ in the sense that $\Sbep[X_n]\searrow 0$ whenever $\mathscr{H}_b\ni X_n\searrow 0$ (cf. Zhang \cite[Lemma 2.5]{ZhangLLN} and Zhang \cite[Lemma 4.2]{ZhangLIL21}). In particular, if $\Omega$ is a Polish space, $\mathscr{H}=C(\Omega)$ is the space of continuous functions on $\Omega$, and $\mathscr{P}$ is a weakly compact family of probabilities on $\Omega$, then $\Sbep$ defined by
$\Sbep[\varphi]=\sup_{P\in\mathscr{P}}E_P[\varphi]$ satisfies condition (CC). 

In these paper, we consider all possible precise limit points of $S_n/n$ under this condition, answer Ter\'an's \cite{Teran14} question (Q1) under Peng's independence, and, under \eqref{eq:moment} or Peng's condition \eqref{eq:Pengcondition}, show that \eqref{eq:sec2.1} and \eqref{eq:sec2.2} hold for some $P\in \mathscr{P}$ when condition (CC) is satisfied with $\mathscr{P}$. We introduce a condition which is weaker than both \eqref{eq:Pengcondition} and \eqref{eq:moment}:
\begin{align}\label{eq:weakmoment}
\vSbep[(|X_1|-c)^+]\to 0 \text{ as } c\to \infty.
\end{align}
Let the core of $\Sbep$ be the set of all probability measures dominated by $\Sbep$,
\begin{align*}
\mathrm{core}(\Sbep)=& \left\{P: P \text{ is a probability measure on }\sigma(\mathscr{H}) \right. \\
& \left.\text{ with }
E_P[X]\le \Sbep[X] \text{ for all } X\in \mathscr{H}_b \right\}.
\end{align*}
The following is our basic theorem.

\begin{theorem}\label{th2}
Suppose that \eqref{eq:weakmoment} and condition (CC) with $\mathscr{P}$ are satisfied.
\begin{description}
\item[\rm (a)] For any sequence of functions $\{\varphi_i(x_1,\ldots,x_i);i\ge 1\}$ for which each $\varphi_i(x_1,\ldots, x_i)$ is a Borel-measurable function on $\mathbb R^i$ and $\varphi_i(x_1,\ldots,x_i)\in [\underline{\mu},\overline{\mu}]$, $i=1,2,\ldots$, there is a probability measure $Q\in \mathrm{core}(\Sbep)$ such that
$\{X_n-\varphi_{n-1}(X_1,\ldots,X_{n-1}); n\ge 2\}$ is a sequence of martingale differences under $Q$, and
\begin{align}\label{eqth2.2}
& V\left(\lim_{n\to \infty}\frac{S_n-\sum_{i=1}^n\varphi_i(X_1,\ldots,X_i)}{n}=0\right)\nonumber\\
=& Q\left(\lim_{n\to \infty}\frac{S_n-\sum_{i=1}^n\varphi_i(X_1,\ldots,X_i)}{n}=0\right)=1,
\end{align}
where $V=\upCapc, \outCapc$ or $\Capc^{\mathscr{P}}$.
\item[\rm (b)] For any sequence of continuous functions $\{\varphi_i(x_1,\ldots,x_i)\in [\underline{\mu},\overline{\mu}]; i\ge 1\}$, there is a probability measure $P\in \mathscr{P}$ such that
\begin{align}\label{eqth2.3}
& V\left(\lim_{n\to \infty}\frac{S_n-\sum_{i=1}^n\varphi_i(X_1,\ldots,X_i)}{n}=0\right)\nonumber\\
=& P\left(\lim_{n\to \infty}\frac{S_n-\sum_{i=1}^n\varphi_i(X_1,\ldots,X_i)}{n}=0\right)=1,
\end{align}
where $V=\upCapc, \outCapc$ or $\Capc^{\mathscr{P}}$.
\item[\rm (c)] Without condition (CC), there exists a sequence $\delta_n\searrow 0$ such that for any sequence of continuous functions $\{\varphi_i(x_1,\ldots,x_i)\in [\underline{\mu},\overline{\mu}]; i\ge 1\}$,
\begin{align}\label{eqth2.6}
\min_{N\ge n} \upCapc\left(\bigcap_{m=n}^N \left\{ \frac{|S_m-\sum_{i=1}^m\varphi_i(X_1,\ldots,X_i)|}{m}\le \delta_m\right\}\right)\to 1,
\end{align}
and moreover,
\begin{align}\label{eqth2.4}
&\min_{N\ge n} \lowCapc\left(\bigcap_{m=n}^N \left\{ \underline{\mu}-\delta_m\le \frac{S_m}{m} \le \overline{\mu}+\delta_m \right\}\right)\to 1, \\
\label{eqth2.5}
&\quad \outcCapc\left( \underline{\mu}-\delta_n\le \frac{S_n}{n} \le \overline{\mu}+\delta_n \text{ eventually} \right)=1
\end{align}
if \eqref{eq:moment} is satisfied.
\end{description}
\end{theorem}

Based on this theorem, we have the following results on the limit points of the law of large numbers.

\begin{theorem}\label{th3}
Suppose that \eqref{eq:weakmoment} and condition (CC) with $\mathscr{P}$ are satisfied.
\begin{description}
\item[\rm (a)] For any finite-dimensional Borel-measurable functions $\underline{\varphi}(x_1,\ldots,x_d)$ and $\overline{\varphi}(x_1,\ldots,x_d)$ with
$\underline{\mu}\le \underline{\varphi}(x_1,\ldots,x_d)\le \overline{\varphi}(x_1,\ldots,x_d)\le \overline{\mu}$, there is a probability measure $P\in \mathrm{core}(\Sbep)$ such that
\begin{align}\label{eqth3.1}
&V\left(C\Big\{\frac{S_n}{n}\Big\}=\Big[\underline{\varphi}(X_1,\ldots,X_d),\overline{\varphi}(X_1,\ldots,X_d) \Big]\right)\nonumber\\
=& P\left(C\Big\{\frac{S_n}{n}\Big\}=\Big[\underline{\varphi}(X_1,\ldots,X_d),\overline{\varphi}(X_1,\ldots,X_d) \Big]\right)=1,
\end{align}
where $C\{x_n\}$ denotes the cluster set of a sequence $\{x_n\}$ in $\mathbb R$, and $V=\upCapc, \outCapc$ or $\Capc^{\mathscr{P}}$.
\item[\rm (b)] For any continuous functions $\underline{\varphi}(x_1,x_2,\ldots)$ and $\overline{\varphi}(x_1,x_2,\ldots)$ on $\mathbb R^{\mathbb N}$ with
$\underline{\mu}\le \underline{\varphi}(x_1,x_2,\ldots)\le \overline{\varphi}(x_1,x_2,\ldots)\le \overline{\mu}$, there is a probability measure $P\in \mathscr{P}$ such that
\begin{align}\label{eqth3.2}
&V\left(C\Big\{\frac{S_n}{n}\Big\}=\Big[\underline{\varphi}(\bm X),\overline{\varphi}(\bm X) \Big]\right)\nonumber\\
=& P\left(C\Big\{\frac{S_n}{n}\Big\}=\Big[\underline{\varphi}(\bm X),\overline{\varphi}(\bm X) \Big]\right)=1,
\end{align}
\end{description}
where $V=\upCapc, \outCapc$ or $\Capc^{\mathscr{P}}$.
\end{theorem}

\eqref{eqth3.1} and \eqref{eqth3.2} tell us that if condition (CC) is satisfied, then \eqref{eq:sec2.1} and \eqref{eq:sec2.2} hold simultaneously with upper capacity $1$, so that \eqref{eq:PengMC} remains true as Peng \cite{peng2009,peng2010} claimed; furthermore, they can hold simultaneously with upper capacity $1$ with $\overline{\mu}$ and $\underline{\mu}$ being replaced by random variables $\overline{\varphi}(\bm X)$ and $\underline{\varphi}(\bm X)$. Thus, $\Sbep[\varphi(X_1)]$ is a limit value but not the unique limit value of
$\hat{\mathbb M}[\varphi]$ in \eqref{eq:PengMC}.

By choosing $\underline{\varphi}(\bm x)\equiv\underline{\mu}$ and $\overline{\varphi}(\bm x)\equiv\overline{\mu}$, and $\underline{\varphi}(\bm x)=\overline{\varphi}(\bm x)\equiv\mu$, respectively, we obtain the following corollary.

\begin{corollary} \label{thLLN1}
Suppose that \eqref{eq:weakmoment} and condition (CC) with $\mathscr{P}$ are satisfied. We have, for $V=\upCapc, \outCapc$ and $\Capc^{\mathscr{P}}$,
\begin{equation}\label{eq:thLLN1.2}
V\left( C\Big\{\frac{S_n}{n}\Big\}=[\underline{\mu},\overline{\mu}]\right)=1,
\end{equation}
and
\begin{equation}\label{eq:thLLN1.3}
V\left(\lim_{n\to\infty}\frac{S_n}{n}=\mu\right)=1 \quad \text{for all } \mu\in [\underline{\mu},\overline{\mu}].
\end{equation}
Furthermore, there exist $P, P_{\mu}\in \mathscr{P}$ such that
\begin{equation}\label{eq:thLLN1.4}
P\left( C\Big\{\frac{S_n}{n}\Big\}=[\underline{\mu},\overline{\mu}]\right)=1,
\end{equation}
\begin{equation}\label{eq:thLLN1.5}
P_{\mu} \left(\lim_{n\to\infty}\frac{S_n}{n}=\mu\right)=1 \quad \text{for all } \mu\in [\underline{\mu},\overline{\mu}].
\end{equation}
\end{corollary}

\begin{remark}\label{rk:1}
The counterexample given by Ter\'an \cite{Teran18} shows that \eqref{eq:thLLN1.2} and \eqref{eq:thLLN1.3} may fail if condition (CC) is not satisfied (cf. Theorem 3.1 of \cite{Teran18}).
Song \cite{Song2022} has also shown \eqref{eq:SLLN0} and \eqref{eq:thLLN1.5} under the moment condition that $\Sbep[|X_1|^{1+\alpha}]<\infty$ when $\Omega$ is a Polish space and $\mathscr{P}$ is a weakly compact family of probability measures on $\Omega$. For the special case that $\Omega=\mathbb R^{\mathbb N}$ and $X_i(\bm x)=x_i$ for $\bm x=(x_1,x_2,\ldots)$, under $\Sbep[|X_1|^{1+\alpha}]<\infty$, Song \cite{Song2022} also proved that the constant $\mu$ can be a finite-dimensional measurable function $\varphi$. As a direct corollary of Theorem \ref{th3}, we have a similar result as follows for a general space $\Omega$.
\end{remark}

\begin{corollary}\label{cor1}
Suppose that \eqref{eq:weakmoment} and condition (CC) with $\mathscr{P}$ are satisfied.
\begin{description}
\item[\rm (a)] For any finite-dimensional Borel-measurable function $\varphi(x_1,\ldots,x_d):\mathbb R^d\to [\underline{\mu},\overline{\mu}]$, there is a probability measure $P_{\varphi}\in \mathrm{core}(\Sbep)$ such that
\begin{align}\label{eqcor1.2}
&V\left(\lim_{n\to \infty}\frac{S_n}{n}= \varphi(X_1,\ldots,X_d) \right)\nonumber\\
=& P_{\varphi}\left(\lim_{n\to \infty}\frac{S_n}{n}= \varphi(X_1,\ldots,X_d) \right)=1,
\end{align}
where $V=\upCapc, \outCapc$ or $\Capc^{\mathscr{P}}$.
\item[\rm (b)] For any continuous function $\varphi (x_1,x_2,\ldots):\mathbb R^{\mathbb N}\to [\underline{\mu},\overline{\mu}]$, there is a probability measure $P_{\varphi}\in \mathscr{P}$ such that
\begin{align}\label{eqcor1.3}
V\left(\lim_{n\to \infty}\frac{S_n}{n}= \varphi(\bm X) \right)
= P_{\varphi}\left(\lim_{n\to \infty}\frac{S_n}{n}= \varphi(\bm X) \right)=1,
\end{align}
where $V=\upCapc, \outCapc$ or $\Capc^{\mathscr{P}}$.
\end{description}
\end{corollary}

\begin{remark}\label{rk:2}
If \eqref{eqcor1.2} or \eqref{eqcor1.3} holds for all $P$ ($\in\mathrm{core}(\Sbep)$ or $\mathscr{P}$) simultaneously, then by Corollary 3.3 of Zhang \cite{ZhangLLN} we must have $\underline{\mu}=\overline{\mu}$ and $P(\varphi(\bm X)=\overline{\mu})$ for all $P$. Thus, \eqref{eqKol1} cannot hold for all $P$, but it holds at least for one $P$ when $b\in [\underline{\mu}, \overline{\mu}]$.
\end{remark}

\begin{remark}\label{rk:3}
In \eqref{eqth3.2} and \eqref{eqcor1.3}, the continuous functions $\underline{\varphi}$, $\overline{\varphi}$, and $\varphi$ can be replaced by functions in $\mathbb L_c$, where $\mathbb L_c$ is the completion of $C_b(\mathbb R^{\mathbb N})$ under $\mathbb E^{\mathscr{P}}$ in the sense that $\varphi\in \mathbb L_c$ if and only if there is $\varphi_n\in C_b(\mathbb R^{\mathbb N})$ for which $\mathbb E^{\mathscr{P}}|\varphi_n(\bm X)-\varphi(\bm X)|\to 0$. By Theorem 6.1.29 of Peng \cite{peng2019},
\begin{align*}
\mathbb L_c=& \Big\{f(\bm x): \mathbb E^{\mathscr{P}}[|f(\bm X)|]<\infty, \lim_{c\to \infty}\mathbb E^{\mathscr{P}}[(|f(\bm X)|-c)^+]=0, \\
& \text{there is a function } g(\bm x) \text{ such that } \Capc^{\mathscr{P}}(g(\bm X)\ne f(\bm X))=0 \text{ and}, \forall \epsilon>0, \\
& \text{there exists an open set } O \text{ with } \Capc^{\mathscr{P}}(\bm X\in O)<\epsilon
\text{ such that } g\big|_{O^c} \text{ is continuous} \Big\}.
\end{align*}
\end{remark}

\begin{remark}\label{rk:4}
If we let $\mathbb L$ be the completion of $\{f: f $  is a Borel-measurable function on $\mathbb R^d$, $d\ge 1\}$ under $\mathbb E^{\mathrm{core}(\Sbep)}$, then in \eqref{eqth3.1} and \eqref{eqcor1.2}, the finite-dimensional measurable functions $\underline{\varphi}$, $\overline{\varphi}$, and $\varphi$ can be replaced by functions in $\mathbb L$.
We conjecture that they can be replaced by Borel-measurable functions $\varphi(x_1,x_2,\ldots):\mathbb R^{\mathbb N}\to [\underline{\mu},\overline{\mu}]$. The problem is that $\mathbb L$ may not coincide with the space of all Borel-measurable functions on $\mathbb R^{\mathbb N}$.
\end{remark}

\begin{proposition} \label{th:example}
There is a weakly compact family of probabilities $\mathscr{P}$ on $\mathbb R^{\mathbb N}$ and a bounded Borel function $f$ on $\mathbb R^{\mathbb N}$ such that $f$ is not a limit of finite-dimensional Borel functions $f_n$ in the sense that $\mathbb E^{\mathscr{P}}|f_n-f|\to 0$. Moreover, the canonical process $\{X_n;n\ge 1\}$ is a sequence of i.i.d. random variables under $\mathbb E^{\mathscr{P}}$, where $X_i(\omega)=\omega_i$ for $\omega=(\omega_1,\omega_2,\ldots)\in \mathbb R^{\mathbb N}$.
\end{proposition}

Next, we consider the weak law of large numbers.

\begin{corollary} \label{cor2}
Let $\{X_n;n\ge 1\}$ be a sequence of i.i.d. random variables in a sub-linear expectation space $(\Omega,\mathscr{H},\Sbep)$ with \eqref{eq:weakmoment}.
\begin{description}
\item[\rm (a)] If condition (CC) is satisfied, then for any finite-dimensional Borel-measurable function $\varphi(x_1,\ldots,x_d):\mathbb R^d\to [\underline{\mu},\overline{\mu}]$, there is a probability measure $P_{\varphi}\in \mathrm{core}(\Sbep)$ such that
\begin{align}\label{eqcor2.1}
& P_{\varphi}\left(\Big|\frac{S_n}{n}- \varphi(X_1,\ldots,X_d)\Big|\le \epsilon \right)\to 1\quad \text{for all } \epsilon>0.
\end{align}
\item[\rm (b)] If condition (CC) is satisfied with $\mathscr{P}$, then for any continuous function $\varphi (x_1,x_2,\ldots):\mathbb R^{\mathbb N}\to [\underline{\mu},\overline{\mu}]$, there is a probability measure $P_{\varphi}\in \mathscr{P}$ such that
\begin{align}\label{eqcor2.2}
P_{\varphi}\left(\Big|\frac{S_n}{n}- \varphi(\bm X)\Big|\le \epsilon \right)\to 1 \quad \text{for all } \epsilon>0.
\end{align}
\item[\rm (c)] Without condition (CC), for a sub-additive capacity $\Capc$ satisfying \eqref{eq1.4} we have
\begin{equation}\label{eqcor2.3}
\cCapc\left(\underline{\mu}-\epsilon \le \frac{S_n}{n}\le \overline{\mu}+\epsilon \right)\to 1 \quad \text{for all } \epsilon>0,
\end{equation}
and
\begin{align}\label{eqcor2.4}
\Capc\left(\Big|\frac{S_n}{n}- \varphi(\bm X)\Big|\le \epsilon \right)\to 1 \quad \text{for all } \epsilon>0
\end{align}
for any continuous function $\varphi(X_1,\ldots,X_d):\mathbb R^d\to [\underline{\mu},\overline{\mu}]$, or uniformly continuous function $\varphi(\bm x):\mathbb R^{\mathbb N}\to [\underline{\mu},\overline{\mu}]$. Furthermore, if $\Capc$ is a countably sub-additive capacity which satisfies \eqref{eq1.4}, then we have \eqref{eqcor2.4} for any continuous function $\varphi(\bm x):\mathbb R^{\mathbb N}\to [\underline{\mu},\overline{\mu}]$.
\end{description}
\end{corollary}

The following corollary tells us that condition (CC) can be replaced by a weaker condition that $\Capc\bm X^{-1}$ is continuous from above on closed sets of $\mathbb R^{\mathbb N}$.

\begin{corollary}\label{rk:5}
Let
\begin{align*}
\mathscr{H}(\bm X) =& \{\xi= \varphi(X_1,\ldots,X_d): \varphi\in C_{b,Lip}(\mathbb R^d)\text{ for some } d\ge 1\},\\
\mathrm{core}(\Sbep|\mathscr{H}(\bm X)) =& \left\{P: P \text{ is a probability measure on }\sigma(\mathscr{H}(\bm X)) \right. \\
& \left.\text{ with }
E_P[Z]\le \Sbep[Z] \text{ for all } Z\in \mathscr{H}(\bm X) \right\}.
\end{align*}
Then Theorems \ref{th2}, \ref{th3} and Corollaries \ref{thLLN1}, \ref{cor1}, \ref{cor2} remain true with $\mathrm{core}(\Sbep|\mathscr{H}(\bm X))$ taking the place of $\mathrm{core}(\Sbep)$ and $\mathscr{P}$, if condition (CC) is replaced by the following weaker condition: there exists a sub-additive capacity $\Capc$ satisfying \eqref{eq1.4} such that $\Capc\bm X^{-1}$ is continuous from above on closed sets of $\mathbb R^{\mathbb N}$ in the sense that
\begin{equation}\label{eq:rk5}
\Capc(\bm X\in C_n)\searrow \Capc(\bm X\in C) \text{ whenever } C_n\searrow C \text{ and } C_n, C \text{ are closed sets of } \mathbb R^{\mathbb N}.
\end{equation}
\end{corollary}

The last theorem shows that when condition (CC) is not satisfied, on the same sub-linear expectation space, the conclusions in Theorem \ref{th3}, Corollary \ref{cor1}, and Corollary \ref{cor2} (a), (b) may hold or may fail.

\begin{theorem}\label{th4}
There exist a sub-linear expectation space $(\Omega,\mathscr{H},\Sbep)$ and a family $\mathscr{P}$ of probability measures on $\Omega$ with \eqref{eq:representbyP} such that, for any $\underline{\mu}<\mu<\overline{\mu}$, there exist two sequences $\{X_n;n\ge 1\}, \{Y_n;n\ge 1\}\subset\mathscr{H}$ satisfying the following properties:
\begin{description}
\item[\rm (a)] Both $\{X_n;n\ge 1\}$ and $\{Y_n;n\ge 1\}$ are sequences of i.i.d. random variables under $\Sbep$, $X_n(\omega), Y_n(\omega)\in\{\underline{\mu},\mu,\overline{\mu}\}$, and $X_n$ and $Y_m$ are identically distributed. Further, for each $n$ and $m$, $\{Z_1,\ldots,Z_{n+m}\}$ are i.i.d. under $\Sbep$, where $\{Z_1,\ldots,Z_{n+m}\}$ is a permutation of $\{X_1,\ldots,X_n,Y_1,\ldots,Y_m\}$.
\item[\rm (b)] $\Sbep[\varphi(X_1)]=\max\{\varphi(\underline{\mu}), \varphi(\mu),\varphi(\overline{\mu})\}$, $\cSbep[\varphi(X_1)]=\min\{\varphi(\underline{\mu}), \varphi(\mu),\varphi(\overline{\mu})\}$, $\varphi\in C_{l,Lip}(\mathbb R)$. In particular,
$\Sbep[X_1]=\overline{\mu}$, $\cSbep[X_1]=\underline{\mu}$.
\item[\rm (c)] For any $g\in C_{l,Lip}(\mathbb R)$,
\begin{equation}\label{eq:th4.0}
\cSbep[g(Y_1)]\le \frac{1}{n}\sum_{i=1}^n g(Y_i(\omega))\le \Sbep[g(Y_1)],\quad \omega\in \Omega,
\end{equation}
\begin{equation}\label{eq:th4.1}
\cSbep[g(X_1)]\le \frac{1}{n}\sum_{i=1}^n g(X_i(\omega))\le \Sbep[g(X_1)],\quad \omega\in \Omega,
\end{equation}
\begin{equation}\label{eq:th4.2}
\frac{1}{n}\sum_{i=1}^n g(X_i(\omega))\to g(\mu),\quad \omega\in \Omega.
\end{equation}
\item[\rm (d)] For any continuous function $\varphi(\bm x):\mathbb R^{\mathbb N}\to [\underline{\mu},\overline{\mu}]$,
\begin{equation}\label{eq:th4.3}
\upCapc\left(\left|\frac{\sum_{i=1}^n X_i}{n}-\varphi(\bm X)\right|\le \epsilon\right)\ge \Capc^{\mathscr{P}}\left(\left|\frac{\sum_{i=1}^n X_i}{n}-\varphi(\bm X)\right|\le \epsilon\right)\to 1 \quad \text{for all } \epsilon>0.
\end{equation}
\item[\rm (e)] For any finite-dimensional Borel-measurable functions $\underline{\varphi}(\bm x),\overline{\varphi}(\bm x):\mathbb R^d\to [\underline{\mu},\overline{\mu}]$ or continuous functions $\underline{\varphi}(\bm x),\overline{\varphi}(\bm x):\mathbb R^{\mathbb N}\to [\underline{\mu},\overline{\mu}]$ with $\underline{\varphi}(\bm x)\le \overline{\varphi}(\bm x)$, there exists $P\in \mathscr{P}$ such that
$$ \upCapc\left(C\Big\{\frac{\sum_{i=1}^n Y_i}{n}\Big\}=\Big[\underline{\varphi}(\bm Y),\overline{\varphi}(\bm Y) \Big]\right)=P\left(C\Big\{\frac{\sum_{i=1}^n Y_i}{n}\Big\}=\Big[\underline{\varphi}(\bm Y),\overline{\varphi}(\bm Y) \Big]\right)=1.$$
\item[\rm (f)] For any continuous function $\varphi(\bm x):\mathbb R^{\mathbb N}\to [\underline{\mu},\overline{\mu}]$, there exists $P\in \mathscr{P}$ such that
$$ \upCapc\left(\left|\frac{\sum_{i=1}^n Y_i}{n}-\varphi(\bm Y)\right|\le \epsilon\right)\ge P\left(\left|\frac{\sum_{i=1}^n Y_i}{n}-\varphi(\bm Y)\right|\le \epsilon\right)\to 1 \quad \text{for all } \epsilon>0.
$$
\item[\rm (g)] There is no family of probability measures $\mathscr{Q}$ such that $\Sbep$ on $\mathscr{H}(\bm X)$ or $\mathscr{H}(\bm Y)$ satisfies condition (CC) with $\mathscr{Q}$, where
$$ \mathscr{H}(\bm Z)=\left\{\varphi(Z_1,\ldots,Z_d):\varphi\in C_{b,Lip}(\mathbb R^d), d\ge 1\right\},\quad \bm Z=\bm X \text{ or } \bm Y. $$
\end{description}
\end{theorem}

In Theorem \ref{th4}, $\{X_n;n\ge 1\}$ and $\{Y_n;n\ge 1\}$ are copies of each other, but their behaviors are quite different. For the random variables $\{X_n;n\ge 1\}$, \eqref{eq:sec2.1} and \eqref{eq:sec2.2} fail, and Peng's claim \eqref{eq:PengMC} fails for any strictly monotone function $\psi$, due to \eqref{eq:th4.2}. Also, for any $P$, \eqref{eqth3.1} and \eqref{eqth3.2} hold if and only if $P(\overline{\varphi}=\underline{\varphi}=\mu)=1$, and \eqref{eqcor1.3} holds if and only if $P(\varphi=\mu)=1$.

\section{Preliminaries and Lemmas}
\setcounter{equation}{0}

To prove the main results in Section \ref{sectMain}, we need some lemmas. We first need some properties of the representation \eqref{eq:representbyP}.

\begin{lemma}\label{lem:0} 
\begin{description}
  \item[\rm (a)] The following three statements are equivalent:
\begin{description}
\item[\rm (i)] Condition (CC) is satisfied with some representation $\mathscr{P}$.
\item[\rm (ii)] Condition (CC) is satisfied with $\mathrm{core}(\Sbep)$.
\item[\rm (iii)] $\Sbep$ is regular in the sense that $\Sbep[X_n]\searrow 0$ whenever $\mathscr{H}_b\ni X_n\searrow 0$.
\end{description}
  \item[\rm (b) ] If condition (CC) is satisfied with the representation $\mathscr{P}$, then $\mathscr{P}\bm Y^{-1}$ is a weakly compact family of probability measures on $\mathbb R^{\mathbb N}$ for any $Y_i\in \mathscr{H}$, $i=1,2,\ldots$.
   \item[\rm (c) ] Let $\Omega$ be a Polish space, $\mathscr{H}$ be the family $C_{Lip}(\Omega)$ of Lipschitz functions on $\Omega$ or the family $C(\Omega)$ of continuous function on $\Omega$. Then condition (CC) is satisfied with $\mathscr{P}$ if and only if \eqref{eq:representbyP} holds with $\mathscr{P}$ being a weakly compact family of probability measures on $\Omega$.      
   \item[\rm (d) ] Let $\Omega$ be a Polish space, $\mathscr{H}$ be $C_{Lip}(\Omega)$   or  $C(\Omega)$. If condition (CC) is satisfied with $\mathscr{P}$, then   
       $\upCapc$, $\outCapc$ and $\Capc^{\mathscr{P}}$ are continuous from above on closed sets in the sense that
       $$ \Capc(C_n)\searrow \Capc(C) \text{ whenever } C_n\searrow C \text{ and } C_n, C \text{ are closed sets }$$
       for $\Capc=\upCapc$, $\outCapc$ and $\Capc^{\mathscr{P}}$.
       
       Conversely, if there exists a pair of capacities $\Capc$ and $\cCapc=1-\Capc$ satisfying \eqref{eq1.4} such that $\Capc$  is continuous from above on closed sets, then condition (CC) is satisfied with some $\mathscr{P}$. 
\end{description}

\end{lemma}

{\bf Proof.} The proof of the equivalence of (i)--(iii) can be found in Lemma 2.5 and Remark 2.6 of Zhang \cite{ZhangLLN}. 

For (b), suppose condition (CC) is satisfied with $\mathscr{P}$, $Y_i\in \mathscr{H}$, $i=1,2,\ldots$, and $P_n\in \mathscr{P}$. We want to show that there exists a $P\in \mathscr{P}$ and a subsequence $\{P_{n_k}\}$ such that $P_{n_k}\bm Y^{-1}\implies P\bm Y^{-1}$. Here and in the sequel, $\implies$ denotes weak convergence of probability measures. Notice that $Y_i^{(l)}\in \mathscr{H}_b$. Let $\bm Z=\{Y_i^{(l)}, i,l=1,2,\ldots\}$. By condition (CC), $\mathscr{P}\bm Z^{-1}$ is weakly compact. Thus, there exist $P\in\mathscr{P}$ and a subsequence $\{P_{n_k}\}$ such that
\begin{align*}
\lim_{k\to \infty} E_{P_{n_k}}\left[\varphi(Y_1^{(l)},\ldots, Y_{m}^{(l)})\right]
= E_{P}\left[\varphi(Y_1^{(l)},\ldots, Y_{m}^{(l)})\right], \;\; \varphi\in C_b(\mathbb R^m), m,l\ge 1.
\end{align*}
Let $\varphi\in C_b(\mathbb R^m)$. Then
\begin{align*}
&\left|E_{P}\left[\varphi(Y_1,\ldots, Y_{m})\right]-E_{P}\left[\varphi(Y_1^{(l)},\ldots, Y_{m}^{(l)})\right]\right|\\
\le & \|\varphi\|_{\infty} \sum_{i=1}^m P\left(|Y_i|>l\right)
\le \|\varphi\|_{\infty} \sum_{i=1}^m E_P\left[\frac{|Y_i|}{l}\wedge 1\right]
\to 0 \text{ as } l\to\infty,
\end{align*}
and
\begin{align*}
&\left|E_{P_{n_k}}\left[\varphi(Y_1,\ldots, Y_{m})\right]-E_{P_{n_k}}\left[\varphi(Y_1^{(l)},\ldots, Y_{m}^{(l)})\right]\right| \\
\le & \|\varphi\|_{\infty} \sum_{i=1}^m E_{P_k}\left[\frac{|Y_i|}{l}\wedge 1\right] \\
= & \|\varphi\|_{\infty} \sum_{i=1}^m E_{P_k}\left[\frac{|Y_i^{(2l)}|}{l}\wedge 1\right]
\to \|\varphi\|_{\infty} \sum_{i=1}^m E_P\left[\frac{|Y_i^{(2l)}|}{l}\wedge 1\right] \text{ as } k\to \infty\\
& \to 0 \text{ as } l\to\infty
\end{align*}
since $\frac{|x|}{l}\wedge 1\in C_b(\mathbb R)$. It follows that
\begin{align*}
\lim_{k\to \infty} E_{P_{n_k}}\left[\varphi(Y_1,\ldots, Y_{m})\right]
= E_{P}\left[\varphi(Y_1,\ldots, Y_{m})\right], \;\; \varphi\in C_b(\mathbb R^m), m\ge 1.
\end{align*}
Hence, $P_{n_k}\bm Y^{-1}\implies P\bm Y^{-1}$.
The proof of (b) is completed. 

The ``if" part of (c) is obvious. For the``only if" part, let $P_n\in \mathscr{P}$. We want to show that there exist a subsequence $\{P_{n_k}\}$ and a $P\in \mathscr{P}$ such that
\begin{align}\label{eq:proof:lem0:4} E_{P_{n_k}}[\varphi]\to E_P[\varphi], \;\; \varphi\in C_b(\Omega). 
\end{align}

We notice the following fact:
\begin{description}\sl
  \item[\rm Fact:] Let $\Omega$ be a  separable metric space. Then there exists a countable subset $\mathscr{F}\subset C(\Omega)$ consisting of bounded Lipschitz functions such that for every compact set $K\subset \Omega$, $\mathscr{F}|_K$ is dense in $C(K)$ with respect to the uniform norm.
\end{description}
Let $\mathscr{F}$ be such a countable set, and $\{K_i\}$ be a sequence of compact sets  with $K_i\nearrow \Omega$. Denote $\bm Y=(Y_1,Y_2,\ldots)$ be a sequence in $\mathscr{H}_b$ 
consisting of all $f^{(2)}(=(-2)\vee f\wedge 2)$, $f\in \mathscr{F}$,    and $f_i(x)=\mathrm{dist}(x,K_i)\wedge 1$, $i=1,2,\ldots$, where $\mathrm{dist}(x,K)=\inf\{\mathrm{dist}(x,y): y\in K\}$.  By condition (CC), there exist a subsequence $\{P_{n_k}\}$ and a $P\in \mathscr{P}$ such that $P_{n_k}\bm Y^{-1}\implies P\bm Y^{-1}$. Then
$$ E_{P_{n_k}}[f^{(2)}]\to E_P[f^{(2)}], f\in  \mathscr{F} \text{ and } 
P_{n_k} f_i^{-1}\implies P f_i^{-1}. $$
Now, suppose $\varphi\in C_b(\Omega)$. Without loss of generality, we can assume $\|\varphi\|_{\infty}<1$. For any given $0<\epsilon<1$, we choose $K_i$ such that $P(K_i^c)<\epsilon/6$, and $f\in \mathscr{F}$ such that
$$\sup_{x\in K_i}|f(x)-\varphi(x)|<\epsilon/4.$$
Since $f-\varphi$ is a continuous function and $K$ is compact, there exists a constant $0<\delta_i<1$ such that
$$|f(x)-\varphi(x)|<\epsilon/2\; \text{ when }\;  \mathrm{dist}(x,K_i)< \delta_i.$$
Then $|f(x)|\le |\varphi(x)|+\epsilon/2<2$ and $f(x)=f^{(2)}(x)$ when $\mathrm{dist}(x,K_i)<\delta_i$. Hence
\begin{align*}
 |f^{(2)}(x)-\varphi(x)| \le & |f^{(2)}(x)-\varphi(x)| I\{\mathrm{dist}(x,K_i)< \delta_i\}+3I\{\mathrm{dist}(x,K_i)\ge \delta_i\}\\
= & |f(x)-\varphi(x)| I\{\mathrm{dist}(x,K_i)< \delta_i\}+3I\{\mathrm{dist}(x,K_i)\ge \delta_i\}\\
<& \epsilon/2+3I\{\mathrm{dist}(x,K_i)\ge \delta_i\}=\epsilon/2+3I\{f_i\ge \delta_i\}.
\end{align*}
It follows that
$$ E_{P}[|f^{(2)}-\varphi|]<\epsilon/2+
3P(f_i\ge \delta_i)\le \epsilon/2+3P(K_i^c)<\epsilon, $$
and 
\begin{align*}
  \limsup_{k\to\infty} & E_{P_{n_k}}[|f^{(2)}-\varphi|]<\epsilon/2+
3\limsup_{k\to \infty} P_{n_k}f_i^{-1}\left(y: y\ge \delta_i\right)\\
\le & \epsilon/2+
3 P f_i^{-1}\left(y: y\ge \delta_i\right)\le \epsilon/2+3P(K_i^c)<\epsilon, 
\end{align*}
since $P_{n_k}f_i^{-1}\implies Pf_i^{-1}$ and $\{y:y\ge \delta_i\}$ is a closed set. 
Hence, 
$$\limsup_{k\to\infty}\left| E_{P_{n_k}}[\varphi]-E_{P}[\varphi]\right|<2\epsilon. $$
\eqref{eq:proof:lem0:4} is proven. 

Next, we show the fact. 
Since $\Omega$ is separable, there is a countable dense subset
$
D=\{x_1,x_2,\dots\}\subset \Omega.
$
For each $n\in\mathbb N$, define
$$
g_n(x)= 1\wedge \mathrm{dist}(x,x_n),\;\; x\in \Omega.
$$
Then each $g_n$ is continuous, $1$-Lipschitz, and satisfies
$
0\le g_n\le 1.$

Let $\mathscr{F}$ be the set of all functions of the form
$$
P(g_{n_1},\dots,g_{n_k}),
$$
where $k\in\mathbb N$, $n_1,\dots,n_k\in\mathbb N$, and $P$ is a polynomial in $k$ variables with rational coefficients. 
Then $\mathscr{F}$ is a countable subset of $C(X)$, and every $f\in\mathscr{F}$ is bounded and Lipschitz.  

Now fix a compact set $K\subset \Omega$. Clearly $\mathscr{F}|_K$ is a $\mathbb Q$-subalgebra of $C(K)$ containing the constant function $1$, where $\mathbb Q$ is the set of rational numbers.  We show that it separates points of $K$. Let $x,y\in K$ with $x\ne y$, and set
$$
r=\mathrm{dist}(x,y)>0.
$$
Since $D$ is dense in $\Omega$, choose $x_n\in D$ such that
$$
\mathrm{dist}(x,x_n)<  \frac r3 \wedge 1.
$$
Then
$$
g_n(x)=\min\{\mathrm{dist}(x,x_n),1\}=\mathrm{dist}(x,x_n)<\frac r3\wedge 1.
$$
On the other hand,
$$
\mathrm{dist}(y,x_n)\ge \mathrm{dist}(x,y)-\mathrm{dist}(x,x_n)>r-\frac r3\wedge 1.
$$
Hence
$$
g_n(y)=\min\{\mathrm{dist}(y,x_n),1\}> g_n(x),
$$
so $g_n(x)\ne g_n(y)$. Therefore $\mathscr{F}|_K$ separates points of $K$.

By the Stone--Weierstrass theorem, a subalgebra of $C(K)$ that contains constant $1$ and separates points is uniformly dense in $C(K)$. Hence
$$
\overline{\mathscr{F}|_K}^{\|\cdot\|_\infty}=C(K).
$$
Thus $\mathscr{F}$ is dense in $C(K)$ for every compact set $K\subset \Omega$. The proof of (c) is completed.

For  (d), suppose condition (CC) is satisfied with $\mathscr{P}$. By (c) of the lemma, $\mathscr{P}$ is a weakly compact family on $\Omega$. By Lemma 6.1.12 of Peng \cite{peng2019}, $\Capc^{\mathscr{P}}$ is continuous from above on closed sets. It is sufficient to show that
\begin{align}\label{eq:proof:lem0.5} \upCapc(C)=\outCapc(C)=\Capc^{\mathscr{P}}(C), \text{ for any closed set } C\subset \Omega. 
\end{align}
$\Capc^{\mathscr{P}}(C)\le \outCapc(C)\le \upCapc(C)$ is obvious. Let $C$ be a closed set. Then $x\in C$ if and only if $d(x,C)=0$. Let $f_k(x)=1- \{k\cdot \mathrm{dist}(x,C)\}\wedge 1$. Then $f_k\in \mathscr{H}_b$, and $I_C\le f_k\searrow I_C$. Let $P_k\in \mathscr{P}$ such that 
$$ E_{P_k}[f_k]> \sup_{P\in \mathscr{P}}E_P[f_k]-\frac{1}{k}. $$
By the compactness of $\mathscr{P}$, there exist a subsequence $\{P_{k_i}\}$ and a $P\in\mathscr{P}$ such that $P_{k_i}\implies P$. 
Then for $i\ge j$,
\begin{align*}
\upCapc(C)\le &\Sbep[f_{k_i}]=\sup_{P\in \mathscr{P}}E_P[[f_{k_i}]
\le E_{P_{k_i}}[f_{k_i}]+\frac{1}{k_i}\\
\le & E_{P_{k_i}}[f_{k_j}]+\frac{1}{k_i}\to E_P[f_{k_j}] \text{ as } i\to \infty
\end{align*}
Letting $j\to \infty$ yields $\upCapc(C)\le P(C)\le \Capc^{\mathscr{P}}(C)$. \eqref{eq:proof:lem0.5} is proven. 

Now, suppose $\Capc$ is continuous from above on closed sets. 
Let $\mathscr{H}_b\ni X_n\searrow 0$. For any $\epsilon>0$, 
$\{X_n\ge \epsilon\}$ is a closed set and $\{X_n\ge \epsilon\}\searrow \emptyset$.
It follows that
$$0\le \Sbep[X_n]\le \epsilon +\Capc(X_n\ge \epsilon)\searrow \epsilon $$
by the assumption. Thus, $\Sbep$ is regular. By (a) of the Lemma, condition (CC) is satisfied. 
The proof is completed.
$\Box$

\begin{lemma}\label{lem:1}
\begin{description}
\item[\rm (i)] Let $\mathscr{P}$ and $\mathscr{Q}$ be two families of probability measures on $\mathbb R^{\mathbb N}$ with
\begin{equation}\label{eq:lem1.1}
\Capc^{\mathscr{P}}(\bm x: |x_i|>c)\to 0 \quad \text{or}\quad \Capc^{\mathscr{Q}}(\bm x: |x_i|>c)\to 0 \text{ as } c\to \infty \text{ for each } i.
\end{equation}
If
\begin{equation}\label{eq:lem1.2}
\sup_{P\in\mathscr{P}}E_P[\varphi]= \sup_{Q\in\mathscr{Q}}E_Q [\varphi] \quad \text{for all } \varphi\in C_{b}^{\infty}(\mathbb R^d), d\ge 1,
\end{equation}
then
\begin{equation}\label{eq:lem1.3}
\sup_{P\in\mathscr{P}}E_P[\varphi]= \sup_{Q\in\mathscr{Q}}E_Q [\varphi] \quad \text{for all } \varphi\in C_{b}(\mathbb R^{\mathbb N}).
\end{equation}
Here and in the sequel, for a $d$-dimensional function $\varphi(x_1,\ldots, x_d)$, we also regard it as an $l$-dimensional function $\varphi^{\ast}(x_1,\ldots,x_l)=\varphi(x_1,\ldots,x_d)$ for $l\ge d$, and as a function $\varphi^{\ast}(x_1,x_2,\ldots)=\varphi(x_1,\ldots,x_d)$ on $\mathbb R^{\mathbb N}$.
\item[\rm (ii)] Let $\mathscr{P}$ and $\mathscr{Q}$ be two convex and weakly compact families of probability measures on $\mathbb R^{\mathbb N}$. If \eqref{eq:lem1.2} holds, then $\mathscr{P}=\mathscr{Q}$.
\item[\rm (iii)] Let $\mathscr{P}$ and $\mathscr{Q}$ be two families of probability measures on $\Omega=C([0,\infty))$, one of which is tight, where $C([0,\infty))$ is the space of continuous functions on $[0,\infty)$ with distance $\mathrm{dist}(x,y)=\sum_{k=1}^{\infty}(1\wedge\|x-y\|_k)/2^k$, $\|x-y\|_k=\sup_{0\le t\le k}|x(t)-y(t)|$. For $t_1,\ldots,t_d \ge 0$, define the projection map $\pi_{t_1,\ldots,t_d}$ from $C([0,\infty))$ to $\mathbb R^d$ by
$\pi_{t_1,\ldots,t_d}x=(x(t_1),\ldots, x(t_d))$. Let $T_0$ be a dense subset of $[0,\infty)$. If
\begin{equation}\label{eq:lem1.4}
\begin{aligned}
& \sup_{P\in\mathscr{P}}E_P[\varphi\circ\pi_{t_1,\ldots,t_d}]= \sup_{Q\in\mathscr{Q}}E_Q [\varphi\circ\pi_{t_1,\ldots,t_d}] \\
& \;\; \text{for all } \varphi\in C_{b}^{\infty}(\mathbb R^d), d\ge 1, t_1,\ldots,t_d\in T_0,
\end{aligned}
\end{equation}
then
\begin{equation}\label{eq:lem1.5}
\sup_{P\in\mathscr{P}}E_P[\varphi]= \sup_{Q\in\mathscr{Q}}E_Q [\varphi] \;\; \text{for all } \varphi\in C_{b}( \Omega),
\end{equation}
and both $\mathscr{P}$ and $\mathscr{Q}$ are tight.
\item[\rm (iv)] Let $\mathscr{P}$ and $\mathscr{Q}$ be two convex and weakly compact families of probability measures on $\Omega=C([0,\infty))$, and let $T_0$ be a dense subset of $[0,\infty)$. If \eqref{eq:lem1.4} holds, then $\mathscr{P}=\mathscr{Q}$.
\end{description}
\end{lemma}

Lemma \ref{lem:1} (i) and (iii) tell us that on the metric spaces $\mathbb R^{\mathbb N}$ and $C([0,\infty))$, a tight sub-linear expectation is determined by its finite-dimensional distributions. (ii) and (iv) tell us that finite-dimensional distributions also determine the convex and weakly compact representation $\mathscr{P}$ of a regular sub-linear expectation on $\mathbb R^{\mathbb N}$ and $C([0,\infty))$.

{\bf Proof.} For (i), it is sufficient to show that, for any $\varphi\in C_b(\mathbb R^{\mathbb N})$ and any $\epsilon>0$, there exist $d$ and a function $\psi_d(x_1,\ldots, x_d)\in C_{b}^{\infty}(\mathbb R^d)$ such that
\begin{equation}\label{eq:prooflem1.1}
\mathbb E^{\mathscr{P}}[|\varphi-\psi_d|]<\epsilon \quad \text{and}\quad \mathbb E^{\mathscr{Q}}[|\varphi-\psi_d|]<\epsilon.
\end{equation}
Let $M=\sup_{\bm x}|\varphi(\bm x)|$. First, by \eqref{eq:lem1.2}, if one of $\mathscr{P}$ and $\mathscr{Q}$ satisfies \eqref{eq:lem1.1}, then they both satisfy \eqref{eq:lem1.1}. Hence, we can choose $c_k$ such that $\Capc^{\mathscr{P}}(|X_k|\ge c_k)< \epsilon /(2^k(4M))$ and $\Capc^{\mathscr{Q}}(|X_k|\ge c_k)< \epsilon /(2^k(4M))$. Let $K=\bigotimes_{k=1}^{\infty}[-c_k,c_k]$. Then $K$ is a compact set in $\mathbb R^{\mathbb N}$. Choose $\bm x_0\in K$. For $\bm x=(x_1,x_2,\ldots)$, denote $\bm x^d=\left((-c_1)\vee x_1\wedge c_1,\ldots,(-c_d)\vee x_d\wedge c_d, x_{d+1}^0, x_{d+2}^0,\ldots \right)$ and $f_d(x_1,\ldots,x_d)=\varphi(\bm x^d)$. Then $f_d$ is a bounded and uniformly continuous function on $\mathbb R^d$, and
$$ \sup_{\bm x}|f_d(\bm x)-\varphi(\bm x)|\le \sup_{\bm x\in K}|\varphi(\bm x^d)-\varphi(\bm x)|+2MI\{\bm x\in K^c\}. $$
Notice that $\varphi(\bm x)$ is uniformly continuous on the compact set $K$, and $\sup_{\bm x }\mathrm{dist}(\bm x^d,\bm x)\le \sum_{j=d+1}^{\infty}1/2^j=1/2^d\to 0$. We have, for $d$ large enough,
$$\sup_{\bm x\in K}|\varphi(\bm x^d)-\varphi(\bm x)| <\epsilon/4.$$
Let $\phi_{d,\delta}(x_1,\ldots,x_d)$ be the density function of the $d$-dimensional normal distribution $N(\bm 0,\delta^2  \bm I)$, and
\begin{align*}
\psi_{\delta,d}(x_1,\ldots,x_d)=&\idotsint \phi_{d,\delta}(x_1-y_1,\ldots,x_d-y_d)f_d(y_1,\ldots,y_d)dy_1\cdots dy_d \\
=&\idotsint \phi_{d,1}(u_1,\ldots,u_d)f_d(x_1-\delta u_1,\ldots,x_d-\delta u_d)d u_1\cdots d u_d
\end{align*}
be the convolution of $\phi_{\delta,d}$ and $f_d$. Then $\psi_{\delta,d}\in C_{b}^{\infty}(\mathbb R^d)$, $|\psi_{\delta,d}(x_1,\ldots,x_d)|\le M$, and
$$\sup_{x_1,\ldots,x_d}|\psi_{\delta,d}(x_1,\ldots,x_d)-f_d(x_1,\ldots,x_d)|\to 0 \; \text{as } \delta\to 0. $$
Choose $\delta>0$ such that $\sup_{x_1,\ldots,x_d}|\psi_{\delta,d}(x_1,\ldots,x_d)-f_d(x_1,\ldots,x_d)|<\epsilon/4$.
Then
\begin{equation}\label{eq:finiteapprox}
\sup_{\bm x}|\varphi(\bm x)-\psi_{\delta,d}(\bm x)|\le \frac{\epsilon}{2}I\{\bm x\in K\}+2MI\{\bm x\in K^c\}.
\end{equation}
It follows that
$$\mathbb E^{\mathscr{P}}[|\psi_{\delta,d}-\varphi|]\le \epsilon/2+2M\sum_{k=1}^{\infty}\Capc^{\mathscr{P}}(\bm x:|x_k|\ge c_k)< \epsilon/2+\epsilon/2\le \epsilon. $$
Similarly, $\mathbb E^{\mathscr{Q}}[|\psi_{\delta,d}-\varphi|]<\epsilon$. \eqref{eq:prooflem1.1} is proved.

For (ii), we first assume that
\begin{equation}\label{eq:prooflem1.2}
\mathbb E^{\mathscr{P}}[(|x_i|-c)^+]\to 0 \text{ as } c\to\infty.
\end{equation}
Notice that for $Q\in \mathscr{Q}$, by \eqref{eq:lem1.2},
\begin{align*}
E_Q[(|x_i|-c)^+]=& \lim_{b\to \infty} E_Q[0\vee (|x_i|-c)\wedge b]\\
\le &\lim_{b\to \infty} \sup_{Q\in \mathscr{Q}}E_Q[0\vee (|x_i|-c)\wedge b]
= \lim_{b\to \infty} \sup_{P\in \mathscr{P}}E_P[0\vee (|x_i|-c)\wedge b]  \\
\le & \sup_{P\in \mathscr{P}}E_P[(|x_i|-c)^+].
\end{align*}
It follows that $\mathbb E^{\mathscr{Q}}[(|x_i|-c)^+]\le \mathbb E^{\mathscr{P}}[(|x_i|-c)^+]$. Thus,
\begin{equation}\label{eq:prooflem1.3}
\mathbb E^{\mathscr{P}}[(|x_i|-c)^+]\to 0 \quad \text{and}\quad \mathbb E^{\mathscr{Q}} [(|x_i|-c)^+]\to 0 \quad \text{as } c\to \infty.
\end{equation}
By (i), \eqref{eq:lem1.3} holds. Notice that for $\varphi\in C_{Lip}(\mathbb R^d)$ we have
$$ \left|\varphi(x_1,\ldots,x_d)-\varphi\left((-c)\vee x_1\wedge c,\ldots,(-c)\vee x_d\wedge c\right)\right|
\le K\sum_{i=1}^d (|x_i|-c)^+. $$
It can be shown that under \eqref{eq:prooflem1.3}, \eqref{eq:lem1.3} implies
\begin{equation}\label{eq:prooflem1.4}
\sup_{P\in\mathscr{P}}E_P[\varphi]= \sup_{Q\in\mathscr{Q}}E_Q [\varphi] \quad \text{for all } \varphi\in C_{Lip}(\mathbb R^d), d\ge 1.
\end{equation}

Suppose $Q\in \mathscr{Q}$. We need to show that $Q\in \mathscr{P}$. By \eqref{eq:prooflem1.4},
$$ E_Q[\varphi]\le \sup_{P\in\mathscr{P}}E_P [\varphi] \quad \text{for all } \varphi\in C_{Lip}(\mathbb R^d). $$
Let $\pi_d(\bm x)=(x_1,\ldots, x_d)$ be the projection map from $\mathbb R^{\mathbb N}$ to $\mathbb R^{d}$. Then $\mathscr{P}\pi_d^{-1}$ is a Wasserstein space $\mathscr{P}_1(\mathbb R^d)$ as in Definition 2.4 of Li and Lin \cite{LiLin2017} by \eqref{eq:prooflem1.3}.
By Corollary 2.8 of Li and Lin \cite{LiLin2017}, $Q|_{\mathbb R^d}\in \mathscr{P}|_{\mathbb R^d}$. That is, $Q\pi_d^{-1}\in \mathscr{P}\pi_d^{-1}$. Hence for each $d$, there exists $P_d\in \mathscr{P}$ such that $Q\pi_d^{-1}=P_d\pi_d^{-1}$. It is obvious that $P_{d^{\prime}}\pi_{d^{\prime}}^{-1}\pi_{d^{\prime},d}^{-1}=Q\pi_{d^{\prime}}^{-1}\pi_{d^{\prime},d}^{-1}=Q\pi_d^{-1}=P_d\pi_d^{-1}$ for $d^{\prime}\ge d$, where $\pi_{d^{\prime},d}(x_1,\ldots,x_{d^{\prime}})=(x_1,\ldots,x_d)$ is the projection map from $\mathbb R^{d^{\prime}}$ to $\mathbb R^d$. Hence, for $d^{\prime}\ge d$,
$$ E_Q[\varphi]=E_{P_d}[\varphi]=E_{P_{d^{\prime}}}[\varphi] \quad \text{for all } \varphi\in C_{b,Lip}(\mathbb R^d). $$
By the weak compactness of $\mathscr{P}$, for the sequence $\{P_{i}\}$ there exist $P\in \mathscr{P}$ and a subsequence $d^{\prime}\to \infty$ such that $P_{d^{\prime}}$ converges weakly to $P$. Hence,
$$ E_Q[\varphi]=\lim_{d^{\prime}\to \infty}E_{P_{d^{\prime}}}[\varphi]=E_P[\varphi] \quad \text{for all } \varphi\in C_{b,Lip}(\mathbb R^d), d\ge 1. $$
It follows that $Q\equiv P\in \mathscr{P}$.

Next, we remove condition \eqref{eq:prooflem1.2}. Notice that $\mathscr{P}$ and $\mathscr{Q}$ are weakly compact and so are tight. Hence \eqref{eq:lem1.1} holds.
By (i), under \eqref{eq:lem1.1}, \eqref{eq:lem1.2} can be extended to
\begin{equation}\label{eq:prooflem1.6}
\sup_{P\in\mathscr{P}}E_P[\varphi]= \sup_{Q\in\mathscr{Q}}E_Q [\varphi] \quad \text{for all } \varphi\in C_b(\mathbb R^d), d\ge 1.
\end{equation}
On the other hand, by \eqref{eq:lem1.1}, for each $i$ there exists a continuous and strictly increasing function $f_i(x)$ such that
$\int_0^{\infty}\Capc^{\mathscr{P}}\big(\bm x: |f_i(x_i)|>x\big)dx<\infty$, which implies $\mathbb E^{\mathscr{P}}[(|f_i(x_i)|-c)^+]\le
\int_c^{\infty}\Capc^{\mathscr{P}}\big(\bm x: |f_i(x_i)|>x\big)dx\to 0$ as $c\to \infty$.
Let $\bm f(\bm x)=(f_1(x_1),f_2(x_2),\ldots)$. Consider two families $\widetilde{\mathscr{P}}=\mathscr{P}\bm f^{-1}$ and $\widetilde{\mathscr{Q}}=\mathscr{Q}\bm f^{-1}$ of probability measures on $\mathbb R^{\mathbb N}$. They are both convex and weakly compact with
\begin{align*}
& \sup_{\widetilde{P}\in \widetilde{\mathscr{P}}}E_{\widetilde{P}}[\varphi(x_1,\ldots,x_d)]\\
=& \sup_{P\in  \mathscr{P}}E_P[\varphi(f_1(x_1),\ldots,f_d(x_d))]=\sup_{Q\in  \mathscr{Q}}E_Q[\varphi(f_1(x_1),\ldots,f_d(x_d))]\\
=&\sup_{\widetilde{Q}\in \widetilde{\mathscr{Q}}}E_{\widetilde{Q}}[\varphi(x_1,\ldots,x_d)]\;\; \text{for all } \varphi\in C_{b,Lip}(\mathbb R^d), \; d\ge 1,
\end{align*}
by \eqref{eq:prooflem1.6}. It follows that $\widetilde{\mathscr{P}}= \widetilde{\mathscr{Q}}$. Notice that $\bm f:\mathbb R^{\mathbb N}\to \mathbb R^{\mathbb N}$ is a one-to-one map. Hence, $\mathscr{P}=\mathscr{Q}$. The proof is completed.

For (iii), we suppose $\mathscr{P}$ is tight.  We show that if $Q\in \mathscr{Q}$, then
\begin{equation}\label{eq:prooflem1.7}
E_Q[\varphi]\le \sup_{P\in\mathscr{P}}E_P [\varphi] \quad \text{for all } \varphi\in C_{b}( \Omega).
\end{equation}

First, we show that for any $\varphi\in C_b(\Omega)$, $\epsilon>0$ and a compact set $K$, there exist  finite many $t_1,t_2,\ldots,t_d\in T_0$ with $0<t_1<t_2<\ldots<t_d$, and $\psi\in C_b(\mathbb R^d)$, such that
\begin{align}\label{eq:finiteapprox}
\left|\psi(\pi_{t_1,\ldots,t_d}x)-\varphi( x)\right|\le \epsilon/2+\|\varphi\|_{\infty}I\{x\not\in K\}. 
\end{align} 

Let $\varphi\in C_b\big(\Omega)$. For $0< t_1<t_2<\ldots<t_d<N$, $t_i\in T_0$, we define a map $\Pi^{-1}_{t_1,\ldots,t_d}$ from $\mathbb R^d$ to $C([0,\infty))$ by
$$\Pi^{-1}_{t_1,\ldots,t_d}(x_1,\ldots,x_d)= \begin{cases}
0,   & \text{if }  t=0 \text{ or } t\ge N; \\
x_k, & \text{if } t=t_k \; (k=1,\ldots, d);\\
\text{linear}, & \text{otherwise}.
\end{cases}$$
Then $\Pi_{t_1,\ldots,t_d}^{-1}$ is a continuous map, $\psi=: \varphi\circ\Pi^{-1}_{t_1,\ldots,t_d}\in C_b(\mathbb R^d)$, and $\|\psi\|_{\infty}\le \|\varphi\|_{\infty}$. Denote $\widetilde{\pi}_{t_1,\ldots,t_d}=\Pi^{-1}_{t_1,\ldots,t_d}\circ \pi_{t_1,\ldots,t_d}$.
  Then $\varphi(\widetilde{\pi}_{t_1,\ldots,t_d} x)=\psi(x(t_1),\ldots, x(t_d))$.

Now, suppose that $t_{i+1}-t_i<\delta$ for $i=0,\ldots, d$, where $t_0=0$, $t_{d+1}=N$. Let $$\omega_{\delta,N}(x)=\sup\limits_{|t-s|<\delta,t,s\in[0,N]}|x(t)-x(s)|.$$ It is easily seen that
$$ \left\|\widetilde{\pi}_{t_1,\ldots,t_d}x -x \right\|_{N-1}\le \left\|\widetilde{\pi}_{t_1,\ldots,t_d}x -x \right\|_{t_d}\le \omega_{\delta,N}(x). $$
Let $\epsilon>0$ be given. Since $\varphi$ is continuous, for each $x$, there is an $\epsilon_x>0$ such that
$$ \left|\varphi(x)-\varphi(y)\right|<\epsilon/4 \text{ whenever } \mathrm{dist}(x,y)<\epsilon_x. $$
Let $K\subset C(\Omega)$ be a compact set. Then it can be covered by a union of finitely many sets $\{y: \mathrm{dist}(x,y)<\epsilon_x/2\}$, $x\in K$. So, there is an $\epsilon_K>0$ such that
$ \left|\varphi(x)-\varphi(y)\right|<\epsilon/2$ whenever $\mathrm{dist}(x,y)<\epsilon_K$ and $x\in K$.  Notice that
\begin{align*}
dist\left(\widetilde{\pi}_{t_1,\ldots,t_d}x, x \right)\le \sum_{k=1}^{N-1}\frac{ \left\|\widetilde{\pi}_{t_1,\ldots,t_d}x -x \right\|_k}{2^k}+\sum_{k=N}^{\infty}1/2^k
\le \omega_{\delta,N}(x)+2/2^N.
\end{align*}
Choose $N$ such that $1/N\le \epsilon_K/4$. It follows that
\begin{align*}
\left|\varphi\left(\widetilde{\pi}_{t_1,\ldots,t_d} x\right)-\varphi(x)\right|
<\frac{\epsilon}{2}I\{x\in K\}+2\|\varphi\|_{\infty} I\{\omega_{\delta,N}(x)\ge \epsilon_K/2, x\in K\}+2\|\varphi\|_{\infty}I\{x\not\in K\}.
\end{align*}
Since $K$ is compact, $\sup_{x\in K}\omega_{\delta,N}(x)\to 0$ as $\delta\to 0$. Hence for $\delta$ small enough, $\{\omega_{\delta,N}(x)\ge \epsilon_K/2, x\in K\}=\emptyset$. Thus
\begin{align*}
\left|\psi(\pi_{t_1,\ldots,t_d}x)-\varphi( x)\right|=\left|\varphi\left(\widetilde{\pi}_{t_1,\ldots,t_d} x\right)-\varphi(x)\right|
<\frac{\epsilon}{2}I\{x\in K\} +2\|\varphi\|_{\infty}I\{x\not\in K\}.
\end{align*} 
\eqref{eq:finiteapprox} is proven. 

Now, 
by the tightness of $\mathscr{P}$ and the probability measure $Q$, we choose the compact set $K$ such that $\Capc^{\mathscr{P}}(K^c)
+Q(K^c)<\epsilon/(4\|\varphi\|_{\infty})$. By \eqref{eq:finiteapprox},
$$ \sup_{P\in \mathscr{P}}E_P\left[\left|\psi\circ\pi_{t_1,\ldots,t_d} -\varphi\right|\right]
\le \epsilon/2+\|\varphi\|_{\infty} \Capc^{\mathscr{P}}(K^c)<\epsilon, $$
and
$$ E_Q\left[\left|\psi\circ\pi_{t_1,\ldots,t_d} -\varphi\right|\right] <\epsilon/2+2\|\varphi\|_{\infty} Q(K^c)<\epsilon. $$
On the other hand,  by \eqref{eq:lem1.4} and the conclusion of (i),
$$ \sup_{P\in \mathscr{P}}E_P[\psi\circ\psi_{t_1,\ldots,t_d}]= \sup_{Q\in \mathscr{Q}}E_Q[\psi\circ\pi_{t_1,\ldots,t_d}]. $$
It follows that
\begin{align*}
E_Q[\varphi]\le & E_Q\left[\psi\circ\pi_{t_1,\ldots,t_d}\right]+\epsilon
\le \sup_{Q\in\mathscr{Q}}E_Q\left[\psi\circ\pi_{t_1,\ldots,t_d}\right]+\epsilon \\
=&\sup_{P\in\mathscr{P}}E_P\left[\psi\circ\pi_{t_1,\ldots,t_d}\right]+\epsilon
\le \sup_{P\in\mathscr{P}}E_P\left[\varphi\right]+2\epsilon.
\end{align*}
\eqref{eq:prooflem1.7} is proved.

 From \eqref{eq:prooflem1.7}, it follows that
\begin{equation}\label{eq:prooflem1.8}
\sup_{Q\in\mathscr{Q}}E_Q[\varphi]\le \sup_{P\in\mathscr{P}}E_P [\varphi] \quad \text{for all } \varphi\in C_{b}( \Omega).
\end{equation}
For a compact set $K$, by \eqref{eq:prooflem1.7},
$$ E_Q\left[\frac{\epsilon\wedge \mathrm{dist}(\cdot,K)}{\epsilon}\right]\le \sup_{P\in\mathscr{P}}E_P\left[\frac{\epsilon\wedge \mathrm{dist}(\cdot,K)}{\epsilon}\right]\le \Capc^{\mathscr{P}}(K^c). $$
Letting $\epsilon\to 0$ yields $Q(K^c)\le \Capc^{\mathscr{P}}(K^c)$. Hence, $\Capc^{\mathscr{Q}}(K^c)=\sup_{Q\in\mathscr{Q}}Q(K^c)\le \Capc^{\mathscr{P}}(K^c)$. It follows that $\mathscr{Q}$ is also tight when $\mathscr{P}$ is tight. Thus, we have
\begin{equation}\label{eq:prooflem1.9}
\sup_{P\in\mathscr{P}}E_P[\varphi]\le \sup_{Q\in\mathscr{Q}}E_Q [\varphi] \quad \text{for all } \varphi\in C_{b}( \Omega).
\end{equation}
\eqref{eq:lem1.5} follows from \eqref{eq:prooflem1.8} and \eqref{eq:prooflem1.9}.

For (iv), we suppose $P\in \mathscr{P}$. We want to show that there exists $Q\in \mathscr{Q}$ such that
$$P(A)=Q(A) \quad \text{for all } A\in  \mathfrak{B}(\Omega), $$
where $\mathfrak{B}(\Omega)$ is the Borel sigma-field on $\Omega$. Without loss of generality, we assume that $T_0$ is countable, and let $T_0=\{t_1,t_2,\ldots\}$. By (ii), there exists $Q\in \mathscr{Q}$ such that
$$ E_P[\varphi\circ \pi_{t_1,\ldots,t_d}]=E_Q[\varphi\circ \pi_{t_1,\ldots,t_d}]\; \text{for all } \varphi\in C_b(\mathbb R^d), d\ge 1. $$
Hence $P$ and $Q$ are equal on the field $\mathfrak{C}=\{\pi_{t_1,\ldots,t_d}^{-1}H: H\in \mathfrak{B}(\mathbb R^d), t_1,\ldots,t_d\in T_0, d\ge 1\}$. Since $T_0$ is dense in $[0,\infty)$, $\sigma(\mathfrak{C})=\mathfrak{B}(\Omega)$. It follows that $P\equiv Q\in \mathscr{Q}$. The proof is completed.
$\Box$

\begin{lemma} \label{lem:2}
Let $(\Omega,\mathscr{H},\Sbep)$ be a sub-linear expectation space and $\{X_n;n\ge 1\}$ be a sequence of random variables on $(\Omega,\mathscr{H},\Sbep)$ with $\upCapc(|X_i|>c)\to 0$ as $c\to \infty$, $i=1,2,\ldots$. Then there exists a convex and weakly compact family $\widetilde{\mathscr{P}}$ of probability measures on the metric space $\widetilde{\Omega}=\mathbb R^{\mathbb N}$ such that
\begin{equation}\label{eq:lem2.0}
\Sbep[\varphi(X_1,\ldots,X_n)]=\sup_{\widetilde{P}\in \widetilde{\mathscr{P}}}E_{\widetilde{P}}[\varphi], \quad \forall \varphi\in C_{b,Lip}(\mathbb R^n), n\ge 1.
\end{equation}
On $\widetilde{\Omega}$, define the canonical process $\{\widetilde{X}_n;n\ge 1\}$ by $\widetilde{X}_i(\widetilde{\omega})=\widetilde{\omega}_i$ for $\widetilde{\omega}=(\widetilde{\omega}_1,\widetilde{\omega}_2,\ldots)$, and the sub-linear expectation $\widetilde{\mathbb E}$ by
$$ \widetilde{\mathbb E}[\varphi]=\sup_{\widetilde{P}\in \widetilde{\mathscr{P}}}E_{\widetilde{P}}[\varphi], \quad \forall \varphi\in C(\mathbb R^{\mathbb N}). $$
Then the sub-linear space $(\widetilde{\Omega},C(\widetilde{\Omega}),\widetilde{\mathbb E})$ satisfies condition (CC) with $\widetilde{\mathscr{P}}$, and $\{\tilde{X}_n;n\ge 1\}$ is a sequence of random variables with $(X_1,\ldots,X_n)\overset{d}=(\tilde X_1,\ldots,\tilde X_n)$ (in the sense of Definition \ref{def1.1} (i)), $n=1,2,\ldots$.

Further, if $X_n\in \{x_1,\ldots,x_d\}$, then $\widetilde{\Omega}$ and $\mathbb R^{\mathbb N}$ can be replaced by $\{x_1,\ldots,x_d\}^{\mathbb N}$.
\end{lemma}

\begin{definition}
The sequence $\{\tilde{X}_n;n\ge 1\}$ in Lemma \ref{lem:2} is called a copy on $\mathbb R^{\mathbb N}$ of $\{X_n;n\ge 1\}$.
\end{definition}

{\bf Proof.} We choose $\widetilde{\mathscr{P}}$ to be
\begin{align}\label{eq:lem2.2}
\widetilde{\mathscr{P}}_{\max}=& \left\{\widetilde{P}: \widetilde{P} \text{ is a probability measure on } \mathbb R^{\mathbb N} \right. \nonumber\\
& \left. \text{ with } E_{\widetilde{P}}[\varphi]\le \Sbep[\varphi\circ \bm X] \text{ for all } \varphi\in \widetilde{\mathscr{H}}_b\right\}
\end{align}
with $\bm X=(X_1,X_2,\ldots)$ and
\begin{equation}\label{eq:lem2.1}
\widetilde{\mathscr{H}}_b=\{\varphi: \varphi\in C_{b,Lip}(\mathbb R^d)\text{ for some } d\ge 1\}.
\end{equation}
Define a sub-linear expectation $\widetilde{\mathbb E}$ on $\widetilde{\mathscr{H}}_b$ by
$$ \widetilde{\mathbb E}[\varphi]=\Sbep[\varphi\circ \bm X] , \quad \varphi\in \widetilde{\mathscr{H}}_b. $$
Then $\widetilde{\mathscr{P}}$ is the core of $\widetilde{\mathbb E}$ on $\widetilde{\mathscr{H}}_b$. Next, we show that $\widetilde{\mathbb E}$ is regular on $\widetilde{\mathscr{H}}_b$.
Suppose $\widetilde{\mathscr{H}}_b\ni Z_n\searrow 0$ and $0\le Z_n\le L$. For each $Z_n$, there exist $d_n$ and $\varphi_n\in C_{b,Lip}(\mathbb R^{d_n})$ such that
$Z_n(\bm x)=\varphi_n(x_1,\ldots,x_{d_n})$. Then $Z_n(\bm x)$ is a continuous function of $\bm x$ on $\mathbb R^{\mathbb N}$. By the condition assumed, for any $\epsilon>0$, there exists $C_i$ such that $\upCapc(|X_i|>C_i)<\epsilon/(2^iL)$. Let $K=\bigotimes_{i=1}^{\infty}[-C_i,C_i]$ and $K_n=\bigotimes_{i=1}^{p_n}[-C_i,C_i]$. Then $K$ is a compact set in $\mathbb R^{\infty}$. Hence
$$ \delta_n=:\sup_{ \bm x\in K}|Z_n(\bm x)|\to 0. $$
Notice that
\begin{align*}
Z_n(\bm x)=\varphi_n(x_1,\ldots,x_{d_n})\le \delta_n+\sup_{\bm x\not \in K}|\varphi_n(x_1,\ldots,x_{d_n})|
\le \delta_n+LI\{(x_1,\ldots,x_{d_n})\not \in K^n\}.
\end{align*}
By the (finite) sub-additivity of $\upCapc$, it follows that
\begin{align*}
\widetilde{\mathbb E}[Z_n]=&\Sbep[\varphi_n(X_1,\ldots,X_{d_n})]\le \delta_n +L\upCapc\left( (X_1,\ldots, X_{d_n})\not \in K^n\right) \\
\le  & \delta_n+L\sum_{i=1}^{d_n}\upCapc(|X_i|>C_i)
\le  \delta_n +\sum_{i=1}^{\infty}\epsilon/2^i\le \delta_n+\epsilon.
\end{align*}
Letting $n\to \infty$ and then $\epsilon\to 0$ yields $\widetilde{\mathbb E}[Z_n]\to 0$. Hence, $\widetilde{\mathbb E}$ is regular on $\widetilde{\mathscr{H}}_b$. Notice that $\mathbb R^{\mathbb N}$ is a countable-dimensional space.
By Lemma \ref{lem:0} (a) and (b), as a core of $\widetilde{\mathbb E}$ on $\widetilde{\mathscr{H}}_b$, $\widetilde{\mathscr{P}}$ is a weakly compact family of probabilities on $\mathbb R^n$. Hence, the sub-linear space $(\widetilde{\Omega},C(\widetilde{\Omega}),\widetilde{\mathbb E})$ satisfies condition (CC) with $\widetilde{\mathscr{P}}$.
The proof is completed. $\Box$

It should be mentioned that a copy $\{\tilde{X}_n;n\ge 1\}$ of $\{X_n;n\ge 1\}$ does not mean that $(X_1,X_2,\ldots)\overset{d}=(\tilde{X}_1,\tilde{X}_2,\ldots)$ because $\Sbep[\varphi(X_1,X_2,\ldots)]$ may have no definition.

\begin{lemma} \label{lem:3}
Let $\Sbep$ satisfy condition (CC), and let $\{X_n;n\ge 1\}$ be a sequence of random variables on $(\Omega,\mathscr{H},\Sbep)$. Let $\bm X=(X_1,X_2,\ldots)$ and let $\widetilde{\mathscr{P}}_{\max}$ be defined as in \eqref{eq:lem2.2}. Then
\begin{equation}\label{eq:identify}
\widetilde{\mathscr{P}}_{\max} =\mathrm{core}(\Sbep)\bm X^{-1}.
\end{equation}
\end{lemma}

{\bf Proof.} This lemma was proved by Guo and Li \cite{GuoLi2023} when $\{X_n;n\ge 1\}$ are i.i.d. with $\Sbep[(|X_1|-c)^+]\to 0$ as $c\to\infty$. Now,
$$ \mathscr{P}\bm X^{-1}\subset \mathrm{core}(\Sbep)\bm X^{-1}\subset \widetilde{\mathscr{P}}_{\max}  $$
implies
\begin{align*}
\Sbep[\varphi\circ\bm X]\ge & \sup_{\widetilde{Q}\in \widetilde{\mathscr{P}}_{\max}}E_{\widetilde{Q}}[\varphi]
\ge \sup_{\widetilde{P}\in \mathrm{core}(\Sbep)\bm X^{-1}}E_{\widetilde{P}}[\varphi]\\
\ge &\sup_{P\in\mathscr{P}}E_P[\varphi\circ\bm X]=\Sbep[\varphi\circ\bm X], \quad \varphi\in C_{b,Lip}(\mathbb R^d).
\end{align*}
Hence
$$ \sup_{\widetilde{Q}\in \widetilde{\mathscr{P}}_{\max}}E_{\widetilde{Q}}[\varphi]
=\sup_{\widetilde{P}\in \mathrm{core}(\Sbep)\bm X^{-1}}E_{\widetilde{P}}[\varphi], \quad \varphi\in C_{b,Lip}(\mathbb R^d).
$$
Notice that $\mathrm{core}(\Sbep)\bm X^{-1}$ is a weakly compact family of probability measures on $\mathbb R^{\mathbb N}$ by Lemma \ref{lem:0} (a) and (b). It is obvious that $\mathrm{core}(\Sbep)\bm X^{-1}$ is convex. By Lemma \ref{lem:2}, $\widetilde{\mathscr{P}}_{\max}$ is also a convex and weakly compact family. By Lemma \ref{lem:1} (ii), \eqref{eq:identify} holds.
$\Box$

\begin{lemma}\label{lem:4}
Suppose that $\Sbep$ satisfies condition (CC), $X\in \mathscr{H}$, and $\vSbep[(|X|-c)^+]\to 0$ as $c\to \infty$. Then
$$ \sup_{P\in \mathscr{P}} E_P[X]=\vSbep[X], $$
and there exists $P\in \mathscr{P}$ such that
$$ E_P[X]=\vSbep[X], $$
$$E_P[\varphi(X)]\le \Sbep[\varphi(X)], \quad \varphi\in C_{b,Lip}(\mathbb R). $$
\end{lemma}

{\bf Proof.} First, for any $c>0$,
$$ \sup_{P\in \mathscr{P}}E_P[X^{(c)}]=\Sbep[X^{(c)}], $$
\begin{align*}
|E_P[X^{(c)}]-E_P[X]|\le & E_P[(|X|-c)^+]=\lim_{b\to \infty}E_P[0\vee(|X|-c)\wedge b]\\
\le & \lim_{b\to \infty}\Sbep[0\vee(|X|-c)\wedge b]=\vSbep[(|X|-c)^+],
\end{align*}
$$ |\vSbep[X]-\Sbep[X^{(c)}]|\le \vSbep[(|X|-c)^+]. $$
Hence
\begin{align*}
& \big|\sup_{P\in \mathscr{P}}E_P[X]-\vSbep[X]\big|\\
\le & |\sup_{P\in \mathscr{P}}E_P[X^{(c)}]-\Sbep[X^{(c)}]|+\sup_{P\in\mathscr{P}}|E_P[X^{(c)}]-E_P[X]|
+|\vSbep[X]-\Sbep[X^{(c)}]|\\
\le & 2\vSbep[(|X|-c)^+]\to 0\text{ as } c\to\infty.
\end{align*}
It follows that $\sup_{P\in \mathscr{P}}E_P[X]=\vSbep[X]$. Choose $P_n\in \mathscr{P}$ such that $|E_{P_n}[X]-\vSbep[X]|\le 1/n$. Notice the weak compactness of $\mathscr{P}$. There exists a subsequence $n^{\prime}$ and $P\in \mathscr{P}$ such that $P_{n^{\prime}}$ converges weakly to $P$. Notice
$$ E_{P_{n^{\prime}}}[X^{(c)}]\to E_P[X^{(c)}] $$
and
\begin{align*}
\big|E_P[X]-\vSbep[X]\big|\le &\big|E_P[X]-E_P[X^{(c)}]\big|+ \big|E_P[X^{(c)}]-E_{P_{n^{\prime}}}[X^{(c)}]\big|\\
& +\big|E_{P_{n^{\prime}}}[X^{(c)}]-E_{P_{n^{\prime}}}[X]\big| +\big|E_{P_{n^{\prime}}}[X]-\vSbep[X]\big| \\
\le & \big|E_P[X^{(c)}]-E_{P_{n^{\prime}}}[X^{(c)}]\big|+2\vSbep[(|X|-c)^+]+1/n^{\prime}\\
& \to 0 \text{ as } n^{\prime}\to \infty \text{ and then } c\to \infty.
\end{align*}
Finally, $P\in \mathscr{P}$ implies
$$ E_P[\varphi(X)]\le \Sbep[\varphi(X)], \quad \varphi\in C_{b,Lip}(\mathbb R). $$
The proof is completed. $\Box$

The next lemma gives the relation between the sub-linear expectation and the conditional expectation under a probability measure; its proof can be found in Guo and Li \cite{GL21} (see also Hu et al. \cite{HLL19} and Gao et al. \cite{GLL21}).

\begin{lemma}\label{lem:5}
Let $\{X_n;n\ge 1\}$ be a sequence of independent random variables in the sub-linear expectation space $(\Omega,\mathscr{H},\Sbep)$ with \eqref{eq:representbyP}. Denote
$$ \mathcal{F}_n=\sigma(X_1,\ldots,X_n) \quad \text{and}\quad \mathcal{F}_0=\{\emptyset,\Omega\}. $$
Then for each $P\in\mathscr{P}$, we have
$$ E_P\left[\varphi(X_n)|\mathcal{F}_{n-1}\right]\le \Sbep[\varphi(X_n)]\quad a.s., \quad \varphi\in C_{b,Lip}(\mathbb R). $$
\end{lemma}

The next lemma will be used in the proof of Theorem \ref{th3}.

\begin{lemma}\label{lem:7}
There exists a sequence $\{\epsilon_i;i\ge 1\}$ such that
\begin{equation}\label{eqlem7.1}
\epsilon_i=1 \text{ or } 0, \quad \liminf_{n\to \infty}\frac{\sum_{i=1}^n \epsilon_i}{n}=0 \quad \text{and}\quad \limsup_{n\to \infty}\frac{\sum_{i=1}^n \epsilon_i}{n}=1.
\end{equation}
\end{lemma}

{\bf Proof.} The lemma may have appeared in the literature. We first give a constructive proof.
Let $0^a$ denote a consecutive string of $a$ zeros, and $1^b$ a consecutive string of $b$ ones. We recursively choose positive integers $a_k, b_k$ such that
$$
a_1=1,\quad b_k>kA_k,\quad a_{k+1}>(k+1)P_k,
$$
where
$$
A_k=\sum_{j=1}^k a_j,\qquad P_k=\sum_{j=1}^k b_j.
$$
Then define the sequence
$$
(\epsilon_1,\epsilon_2,\dots)
=
0^{a_1}\,1^{b_1}\,0^{a_2}\,1^{b_2}\,0^{a_3}\,1^{b_3}\cdots
$$
That is, first place $a_1$ zeros, then $b_1$ ones, then $a_2$ zeros, then $b_2$ ones, and so on.
Let
$
S_n=\sum_{i=1}^n \epsilon_i.
$

Consider the position at the end of the \(k\)-th block of zeros:
$
n_k=A_k+P_{k-1}.
$
At this point, the sequence contains $P_{k-1}$ ones in total, so
$$
\frac{S_{n_k}}{n_k}
=
\frac{P_{k-1}}{A_k+P_{k-1}}
\le
\frac{P_{k-1}}{a_k}
<
\frac1k.
$$
Thus along these $n_k$,
$$
\frac{S_{n_k}}{n_k}\to 0.
$$
Since \(\frac{S_n}{n}\ge 0\), we have
$$
\liminf_{n\to\infty}\frac{S_n}{n}=0.
$$

Consider the position at the end of the $k$-th block of ones:
$
m_k=A_k+P_k.
$
At this point there are $P_k$ ones in total, while the total number of zeros is $A_k$, so
$$
\frac{S_{m_k}}{m_k}
=
\frac{P_k}{A_k+P_k}
=
1-\frac{A_k}{A_k+P_k}
\ge
1-\frac{A_k}{b_k}
>
1-\frac1k.
$$
Hence along these $m_k$,
$$
\frac{S_{m_k}}{m_k}\to 1.
$$
Also, since $\frac{S_n}{n}\le 1$, we get
$$
\limsup_{n\to\infty}\frac{S_n}{n}=1.
$$
Therefore, there indeed exists a sequence \(\{\epsilon_i\}\) such that \eqref{eqlem7.1} holds.
This completes the proof.

Next, we give a proof by the law of large numbers under the sub-linear expectation. Let $\{Y_n;n\ge 1\}$ be a sequence of i.i.d. normal $N(0,[1,2])$ random variables in a sub-linear expectation space $(\widetilde{\Omega},\widetilde{\mathscr{H}}, \widetilde{\mathbb E})$. Without loss of generality, we can assume that
$\widetilde{\mathbb E}$ satisfies condition (CC) with $\widetilde{\mathscr{P}}$, for otherwise we can consider the copy of $\{Y_n;n\ge 1\}$ on $\mathbb R^{\mathbb N}$ as in Lemma \ref{lem:2}.
Notice $\widetilde{\mathcal{E}}[Y_1^2]=1$, $\widetilde{\mathbb E}[Y_1^2]=2$, and $\widetilde{\mathbb E}[Y_1^4]=3\times 2^2<\infty$. By \eqref{eq:thLLN1.2}, which is proved by Zhang \cite{ZhangLLN}, there exists a $P\in \widetilde{\mathscr{P}}$ such that
$$ P\left(\liminf_{n\to \infty}\frac{\sum_{i=1}^n Y_i^2}{n}=1 \text{ and } \limsup_{n\to \infty}\frac{\sum_{i=1}^n Y_i^2}{n}=2\right)\ge 1/2. $$
On the other hand, $\{Y_i^2-E_P[Y_i^2|Y_1,\ldots,Y_{i-1}];i\ge 1\}$ is a sequence of martingale differences with $E_P[(Y_i^2-E_P[Y_i^2|Y_1,\ldots,Y_{i-1}])^2]\le E_P [Y_i^4]\le \widetilde{\mathbb E}[Y_i^4]=12<\infty$. By the law of large numbers for martingales (cf. Theorem 2.18 of Hall and Heyde \cite{HallHeyde80}),
$$ P\left(\lim_{n\to \infty}\frac{\sum_{i=1}^n (Y_i^2-E_P[Y_i^2|Y_1,\ldots,Y_{i-1}])}{n}=0\right)=1. $$
Hence
\begin{align*}
P\Bigg(& \liminf_{n\to \infty}\frac{\sum_{i=1}^n E_P[Y_i^2|Y_1,\ldots,Y_{i-1}]}{n}=1 \\
& \text{and } \limsup_{n\to \infty}\frac{\sum_{i=1}^n E_P[Y_i^2|Y_1,\ldots,Y_{i-1}]}{n}=2\Bigg)\ge 1/2.
\end{align*}
Notice $1=\widetilde{\mathcal E}[Y_i^2]\le E_P[Y_i^2|Y_1,\ldots,Y_{i-1}]\le \widetilde{\mathbb E}[Y_i^2]=2$ a.s. under $P$ by Lemma \ref{lem:5}. Thus, there exists a sequence of real numbers
$\sigma_i\in[1,2]$ such that
\begin{equation}\label{eqprooflem5.2}
\liminf_{n\to \infty} \frac{\sum_{i=1}^n \sigma_i}{n}=1 \quad \text{and}\quad \limsup_{n\to \infty} \frac{\sum_{i=1}^n \sigma_i}{n}=2.
\end{equation}
Let $\epsilon_i=1$ if $\sigma_i>3/2$ and $=0$ if $\sigma_i\le 3/2$.
Notice
$$1\le \frac{1}{2}\epsilon_i+1\le \sigma_i\le \frac{3}{2}+\frac{1}{2}\epsilon_i\le 2. $$
\eqref{eqprooflem5.2} implies \eqref{eqlem7.1}. The proof is completed. $\Box$

\section{Proofs}\label{sec:proof}
\setcounter{equation}{0}

In this section, we give the proofs of the main results in Section \ref{sectMain}.

{\bf Proof of Theorem \ref{th2}.} Let $\mathcal F_i=\sigma(X_1,\ldots, X_i)$, $\mathcal F_0=\{\Omega,\emptyset\}$.

{\em Step 1.} Let $ Q^{\ast}\in \mathrm{core}(\Sbep)$ and
\begin{align*}
\mathscr{Q}(Q^{\ast})=\bigg\{  Q\in \mathrm{core}(\Sbep): Q(X_i\in A)\le 2 Q^{\ast}(X_1\in A), \; A\in \mathscr{B}^1, i\ge 2\bigg\}.
\end{align*}
We first show that for all $Q\in \mathscr{Q}(Q^{\ast})$,
\begin{equation}\label{eqproofth2.2}
Q\left(\lim_{n\to \infty} \frac{\sum_{i=1}^n(X_i-E_Q[X_i|\mathcal F_{i-1}])}{n} =0\right)=1.
\end{equation}
Let $Y_i=X_i^{(i)}=(-i)\vee X_i\wedge i$. Then $\{Y_i-E_Q[Y_i|\mathcal F_{i-1}]; i\ge 1\}$ is a sequence of martingale differences under $Q$ and
\begin{align*}
&\sum_{i=3}^{\infty} \frac{E_Q[(Y_i-E_Q[Y_i|\mathcal F_{i-1}])^2]}{i^2}\le \sum_{i=3}^{\infty}\frac{E_Q[Y_i^2]}{i^2}\le 2\sum_{i=3}^{\infty}\frac{E_{Q^{\ast}}[X_1^2\wedge i^2]}{i^2}\\
=& 2 E_{Q^{\ast}}\bigg[\sum_{i=3}^{\infty}\frac{X_1^2\wedge i^2}{i^2}\bigg]
\le 6 E_{Q^{\ast}}[|X_1|]\le 6\vSbep[|X_1|]<\infty.
\end{align*}
By the law of large numbers for martingales (cf. Theorem 2.18 of Hall and Heyde \cite{HallHeyde80}),
\begin{equation}
Q\left(\lim_{n\to \infty} \frac{\sum_{i=1}^n(Y_i-E_Q[Y_i|\mathcal F_{i-1}])}{n} =0\right)=1.
\end{equation}
On the other hand, by Lemma \ref{lem:5},
\begin{align*}
&\left|E_Q[X_i|\mathcal F_{i-1}]-E_Q[Y_i|\mathcal F_{i-1}]\right|\le E_Q\left[(|X_i|-i)^+|\mathcal F_{i-1}\right]\\
=&  \lim_{b\to \infty}E_Q\left[0\vee (|X_i|-i)\wedge b|\mathcal F_{i-1}\right]\le \lim_{b\to \infty}\Sbep\left[0\vee (|X_i|-i)\wedge b\right]\;\; a.s. \\
=& \lim_{b\to \infty}\Sbep\left[0\vee (|X_1|-i)\wedge b\right]=\vSbep[(|X_1|-i)^+] \to 0 \text{ as } i\to \infty.
\end{align*}
Hence
$$ Q\left(\lim_{n\to \infty} \frac{\sum_{i=1}^n (E_Q[X_i|\mathcal F_{i-1}]-E_Q[Y_i|\mathcal F_{i-1}])}{n}=0\right)=1. $$
Also
\begin{align*}
\sum_{i=2}^{\infty}Q(Y_i\ne X_i)\le \sum_{i=2}^{\infty} 2Q^{\ast}(|X_i|>i)\le 2E_{Q^{\ast}}[|X_1|]\le 2\vSbep[|X_1|]<\infty,
\end{align*}
which implies $Q(Y_i\ne X_i, i.o.)=0$. \eqref{eqproofth2.2} is proved.

Next, we show a more precise result as follows. For any $ Q^{\ast}\in \mathrm{core}(\Sbep)$, there exists a sequence $\epsilon_k\searrow 0$ such that
\begin{equation}\label{eqconvergence}
M=:\sum_{k=1}^{\infty} \sup_{Q\in \mathscr{Q}(Q^{\ast})} Q\left(\max_{ n\le n_{k+1}}\left|\sum_{i=1}^n(X_i-E_Q[X_i|\mathcal F_{i-1}])\right|\ge \epsilon_k n_k\right)<\infty,
\end{equation}
where $n_k=2^k$.

It is sufficient to show that
\begin{equation}\label{eqconvergence1}
\sum_{k=1}^{\infty} \sup_{Q\in \mathscr{Q}(Q^{\ast})} Q\left(\max_{ n\le n_{k+1}}\left|\sum_{i=2}^n(X_i-E_Q[X_i|\mathcal F_{i-1}])\right|\ge \epsilon_k n_k\right)<\infty.
\end{equation}

It is obvious that
\begin{align}\label{eqconvergenceQ}
&\sum_{k=1}^{\infty} \sup_{Q\in \mathscr{Q}(Q^{\ast})} Q\left(\max_{ n\le n_{k+1}}\left|\sum_{i=2}^n(Y_i-E_Q[Y_i|\mathcal F_{i-1}])\right|\ge \epsilon n_k\right) \nonumber\\
\le & \sum_{k=1}^{\infty} \sup_{Q\in \mathscr{Q}(Q^{\ast})} \frac{\sum_{i=2}^{n_{k+1}}E_Q[Y_i^2]}{(\epsilon n_k)^2}\le 2 \sum_{k=1}^{\infty}\sum_{i=2}^{n_{k+1}} \frac{E_{Q^{\ast}}[X_1^2\wedge i^2]}{\epsilon^2 n_k^2}
\nonumber\\
\le & 2\frac{1}{\epsilon^2} \sum_{i=2}^{\infty} E_{Q^{\ast}}[X_1^2\wedge i^2]\sum_{k: n_{k+1}\ge i}n_k^{-2} \le \frac{32}{\epsilon^2} \sum_{i=2}^{\infty} \frac{ E_{Q^{\ast}}[X_1^2\wedge i^2]}{i^2}<\infty,
\end{align}
cf. Corollary 2.1 of Hall and Heyde \cite{HallHeyde80}. Notice
$|E_Q[X_i|\mathcal F_{i-1}]-E_Q[Y_i|\mathcal F_{i-1}]|\overset{a.s.}\le \vSbep[(|X_1|-i)^+]\to 0.$
We have
\begin{equation}\label{eqconvergenceQ1}
\sum_{k=1}^{\infty} \sup_{Q\in \mathscr{Q}(Q^{\ast})} Q\left(\max_{ n\le n_{k+1}}\left|\sum_{i=2}^n(Y_i-E_Q[X_i|\mathcal F_{i-1}])\right|\ge \epsilon n_k\right)<\infty.
\end{equation}
On the other hand,
$$\max_{n_k+1\le n\le n_{k+1}}\left|\sum_{i=2}^n(X_i-Y_i)\right|
\le \sum_{i=1}^{n_{k+1}}(|X_i|-i)^+  $$
and
$$E_Q\left[(|X_i|-i)^+|\mathcal{F}_{i-1}\right]\overset{a.s.}\le \vSbep[(|X_i|-i)^+]\to 0. $$
Write $Z_i=(|X_i|-i)^+$. For $k$ large enough,
\begin{align*}
&  Q\left(\max_{ n\le n_{k+1}}\left|\sum_{i=2}^n(X_i-Y_i)\right|\ge \epsilon n_k\right) \\
\le &  Q\left(\sum_{i=2}^{n_{k+1}}(|X_i|-i)^+\ge \epsilon n_k\right)
\le Q\left(\sum_{i=2}^{n_{k+1}}(Z_i\wedge n_k)\ge \epsilon n_k\right)\\
\le &
Q\left(\sum_{i=2}^{n_{k+1}}\big\{(Z_i\wedge n_k)-E_Q[(Z_i\wedge n_k)|\mathcal{F}_{i-1}]\big\}\ge \epsilon n_k/2\right)\\
\le & \frac{4}{\epsilon^2} \frac{\sum_{i=2}^{n_{k+1}} E_Q[(Z_i\wedge n_k)^2]}{n_k^2}
\le \frac{16}{\epsilon^2}\frac{E_{Q^{\ast}}[X_1^2\wedge n_k^2]}{n_k}\le \frac{32}{\epsilon^2}\sum_{i=n_k+1}^{n_{k+1}}\frac{E_{Q^{\ast}}[X_1^2\wedge i^2]}{i^2}.
\end{align*}
It follows that
$$\sum_{k=1}^{\infty} \sup_{Q\in \mathscr{Q}(Q^{\ast})}  Q\left(\max_{ n\le n_{k+1}}\left|\sum_{i=2}^n(X_i-Y_i)\right|\ge \epsilon n_k\right)<\infty,  $$
which, together with \eqref{eqconvergenceQ1}, yields that
\begin{equation}\label{eqconvergenceQ2}
\sum_{k=1}^{\infty} \sup_{Q\in \mathscr{Q}(Q^{\ast})} Q\left(\max_{ n\le n_{k+1}}\left|\sum_{i=2}^n(X_i-E_Q[X_i|\mathcal F_{i-1}])\right|\ge \epsilon n_k\right)<\infty \text{ for all }\epsilon>0.
\end{equation}
Thus, there exists $\epsilon_k\searrow 0$ such that \eqref{eqconvergence1} holds.

{\em Step 2.} We show that there exists a probability measure $Q^{\ast}\in \mathrm{core}(\Sbep)$, and, for any sequence $\{\varphi_i(x_1,\ldots,x_i), i\ge 1\}$ of functions for which $\varphi_i(x_1,\ldots,x_i)$ is a Borel-measurable function on $\mathbb R^i$ and $\varphi_i(x_1,\ldots,x_i)\in [\underline{\mu},\overline{\mu}]$, there exists a probability measure $Q\in \mathrm{core}(\Sbep)$ such that
\begin{align}\label{eqconditionE}
& E_Q[X_i|\mathcal F_{i-1}]=\varphi_{i-1}(X_1,\ldots,X_{i-1})\;\; a.s. \text{ under } Q, i\ge 2,\\
\label{eq:controlprob}
& Q(X_i\in B)\le 2 Q^{\ast}(X_1\in B), \;\; B\in \mathscr{B}^1, i\ge 2.
\end{align}
Here, the choice of $Q^{\ast}$ does not depend on the $\varphi_i$'s.
Let $\bm X=(X_1,X_2,\ldots)$, and let
$\widetilde{\mathscr{H}}_b$ and $\widetilde{\mathscr{P}}_{\max}$ be defined as in Lemma \ref{lem:3}. By Lemma \ref{lem:3},
$ \widetilde{\mathscr{P}}_{\max} =\mathrm{core}(\Sbep)\bm X^{-1}.
$
Next, we construct $\widetilde{Q}, \widetilde{Q}^{\ast}\in \widetilde{\mathscr{P}}_{\max}$ such that
\begin{align}\label{eqconditionE2}
& E_{\widetilde{Q}}[x_d|x_1,\ldots, x_{d-1}]=\varphi_{d-1}(x_1,\ldots,x_{d-1}), \\
\label{eq:controlprob2}
& \widetilde{Q}(\{\bm x:x_i\in B\})\le 2 \widetilde{Q}^{\ast}(\{\bm x:x_1\in B\}), \;\; B\in \mathscr{B}^1.
\end{align}
Then there will exist $Q, Q^{\ast}\in \mathrm{core}(\Sbep)$ such that $\widetilde{Q}=Q\bm X^{-1}$
and $\widetilde{Q}^{\ast}=Q^{\ast}\bm X^{-1}$ so that \eqref{eqconditionE} and \eqref{eq:controlprob} hold.

For constructing $\widetilde{Q}$, first by Lemma \ref{lem:4} there exist two probability measures $\overline{P}$ and $\underline{P}$ on $\mathbb R$ such that
\begin{equation}\label{eqproofth2.11}
E_{\overline{P}}[\varphi]\le \Sbep[\varphi(X_1)], \quad E_{\underline{P}}[\varphi]\le \Sbep[\varphi(X_1)], \quad \varphi\in C_{b,Lip}(\mathbb R),
\end{equation}
\begin{equation}\label{eqproofth2.12}
E_{\overline{P}}[x]=\vSbep[X_1]=\overline{\mu}, \quad E_{\underline{P}}[x]=\vcSbep[X]=\underline{\mu}.
\end{equation}
For the function $\varphi_{i-1}(x_1,\ldots,x_{i-1})\in [\underline{\mu},\overline{\mu}]$, we let
$$ \alpha_{i-1}=\begin{cases}
\frac{\varphi_{i-1}(x_1,\ldots,x_{i-1})-\underline{\mu}}{\overline{\mu}-\underline{\mu}}, & \text{if } \overline{\mu}\ne \underline{\mu},\\
1, & \text{otherwise}.
\end{cases}
$$
Then $0\le \alpha_{i-1}\le 1$. Define
$$ k_{i-1}(x_1,\ldots,x_{i-1};A)=\alpha_{i-1}\overline{P}(A)+(1-\alpha_{i-1})\underline{P}(A), \;\; A\in \mathscr{B}^1,$$
$$\widetilde{Q}_1 \text{ to be a probability measure on } \mathbb R \text{ satisfying }E_{\widetilde{Q}_1}[\varphi]\le \Sbep[\varphi(X_1)], \varphi\in C_{b,Lip}(\mathbb R), $$
$$\widetilde{Q}_i(A)=\idotsint_{(x_1,\ldots,x_i)\in A} k_{i-1}(x_1,\ldots,x_{i-1}, d x_i)\widetilde{Q}_{i-1}(dx_1,\ldots,d x_{i-1}), \;\; A\in \mathscr{B}^i, i\ge 2,$$
where $\mathscr{B}^i$ is the Borel sigma-field on $\mathbb R^i$.
Then $\{\widetilde{Q}_i; i\ge 1\}$ is a sequence of probability measures which is consistent in the sense that
$$ \widetilde{Q}_i(A\times \mathbb R)=\widetilde{Q}_{i-1}(A), \;\; A\in \mathscr{B}^{i-1}. $$
Hence, by Kolmogorov's extension theorem, there exists a probability measure $\widetilde{Q}$ on $\mathbb R^{\mathbb N}$ such that
$$ \widetilde{Q}\pi_i^{-1}=\widetilde{Q}_i. $$
Now, for any $f(x_1,\ldots,x_d)\in C_{b,Lip}(\mathbb R^d)$,
\begin{align*}
& f_{d-1}(x_1,\ldots,x_{d-1})=: \int f(x_1,\ldots,x_d)k_{d-1}(x_1,\ldots,x_{d-1}; dx_d)\\
=& \alpha_{i-1}\int f(x_1,\ldots,x_d)\overline{P}(dx_d)+(1-\alpha_{i-1})\int f(x_1,\ldots,x_d)\underline{P}(dx_d)\\
\le &\Sbep[f(x_1,\ldots,x_{d-1}, X_d)],
\end{align*}
by \eqref{eqproofth2.11}. It is easily shown that $\Sbep[f(x_1,\ldots,x_{d-1}, X_d)]\in C_{b,Lip}(\mathbb R^{d-1})$. By induction,
\begin{align*}
&E_{\widetilde{Q}}[f(x_1,\ldots,x_d)]=E_{\widetilde{Q}_d}[f(x_1,\ldots,x_d)]=E_{\widetilde{Q}_{d-1}}[f_{d-1}(x_1,\ldots,x_{d-1})]\\
\le & E_{\widetilde{Q}_{d-1}}\left[\Sbep[f(x_1,\ldots,x_{d-1}, X_d)]\right]\le \Sbep\big[\Sbep[f(x_1,\ldots,x_{d-1}, X_d)]\big|_{x_1=X_1,\ldots, x_{d-1}=X_{d-1}}\big]\\
=& \Sbep[f(X_1,\ldots, X_d)].
\end{align*}
It follows that $\widetilde{Q}\in \widetilde{\mathscr{P}}_{\max}$.

On the other hand, for any bounded Borel-measurable function $f(x_1,\ldots,x_{d-1})$,
\begin{align*}
& E_{\widetilde{Q}}[x_df(x_1,\ldots,x_{d-1})]=E_{\widetilde{Q}_d}[x_df(x_1,\ldots,x_{d-1})] \\
=& \idotsint f(x_1,\ldots,x_{d-1})x_dk_{d-1}(x_1,\ldots,x_{d-1};dx_d)\widetilde{Q}_{d-1}(dx_1,\ldots,dx_{d-1})\\
=& \idotsint f(x_1,\ldots,x_{d-1})\big(\alpha_{i-1}\overline{\mu}+(1-\alpha_{i-1})\underline{\mu}\big)\widetilde{Q}_{d-1}(dx_1,\ldots,dx_{d-1})
\;\; (\text{by \eqref{eqproofth2.12}}) \\
=& \idotsint f(x_1,\ldots,x_{d-1})\varphi_{d-1}(x_1,\ldots,x_{d-1})\widetilde{Q}_{d-1}(dx_1,\ldots,dx_{d-1})\\
=&E_{\widetilde{Q}}\left[f(x_1,\ldots,x_{d-1})\varphi_{d-1}(x_1,\ldots,x_{d-1})\right].
\end{align*}
It follows that \eqref{eqconditionE2} holds.

We define $\widetilde{Q}^{\ast}$ by taking $\alpha_i=1/2$ and $\widetilde{Q}_1=\frac{1}{2}\overline{P}+\frac{1}{2}\underline{P}$ in the definition of $\widetilde{Q}$, i.e., $\widetilde{Q}^{\ast}=\big(\frac{1}{2}\overline{P}+\frac{1}{2}\underline{P}\big)^{\bigotimes\mathbb N}$ is the product of $\frac{1}{2}\overline{P}+\frac{1}{2}\underline{P}$. Then $\widetilde{Q}^{\ast}\in \widetilde{\mathscr{P}}_{\max}$ and it does not depend on the $\varphi_i$'s. The canonical process $\{\widetilde{X}_n;n\ge 1\}$ is a sequence of i.i.d. random variables under $\widetilde{Q}^{\ast}$, where $\widetilde{X}_i(\bm x)=x_i$ for $\bm x=(x_1,x_2,\ldots)\in \mathbb R^{\mathbb N}$. Notice that
$$ k_{i-1}(x_1,\ldots,x_{i-1};A)\le 2\big(\frac{1}{2}\overline{P}(A)+\frac{1}{2}\underline{P}(A)\big), \;\; A\in \mathscr{B}^1.$$
\eqref{eq:controlprob2} is obvious.

{\em Step 3.} We show \eqref{eqth2.2} and \eqref{eqth2.3}.

\eqref{eqth2.2} follows from \eqref{eqproofth2.2} and the conclusion in Step 2 immediately. For \eqref{eqth2.3}, let $\varphi_i(x_1,\ldots,x_i)$ be a continuous function with $\varphi_i(x_1,\ldots,x_i)\in [\underline{\mu},\overline{\mu}]$. By the conclusion in Step 2, there exist $Q^{\ast}\in \mathrm{core}(\Sbep)$ and $Q\in \mathscr{Q}(Q^{\ast})$ such that \eqref{eqconditionE} holds. By \eqref{eqconvergence}, there exists a constant $M$ which does not depend on the $\varphi_i$'s such that
\begin{equation}\label{eqconvergence2}
\sum_{k=1}^{\infty} Q\left(\max_{n\le n_{k+1}}\left|\sum_{i=1}^n(X_i-\varphi_{i-1}(X_1,\ldots,X_{i-1}))\right|\ge \epsilon_k n_k\right)\le M<\infty.
\end{equation}
Let $f\in C_{b,Lip}(\mathbb R)$ such that $I\{|x|\ge 2\}\le f(x)\le I\{|x|\ge 1\}$,
$$f_k(\bm x)=f\left(\max_{ n\le n_{k+1}}\left|\sum_{i=1}^n(x_i-\varphi_{i-1}(x_1,\ldots,x_{i-1}))\right|/\epsilon_k n_k\right). $$
Then
\begin{align*}
& E_{Q}\left[\sum_{k=1}^{\infty}f_k(\bm X)\right]\\
\le & \sum_{k=1}^{\infty} Q\left(\max_{ n\le n_{k+1}}\left|\sum_{i=1}^n(X_i-\varphi_{i-1}(X_1,\ldots,X_{i-1}))\right|\ge \epsilon_k n_k\right)\le M<\infty.
\end{align*}
Notice
\begin{align*}
\sup_{P\in \mathrm{core}(\Sbep)}E_P[\varphi(X_1,\ldots,X_d)]=&\sup_{P\in \mathscr{P}}E_P[\varphi(X_1,\ldots,X_d)]\\
=& \Sbep[\varphi(X_1,\ldots,X_d)], \varphi\in C_{b,Lip}(\mathbb R^d).
\end{align*}
Since $\upCapc(|X_i|\ge c)\to 0$ as $c\to\infty$, by Lemma \ref{lem:1} (i) the above inequality can be extended to
$$ \sup_{P\in \mathrm{core}(\Sbep)}E_P[\varphi(X_1,\ldots,X_d)]=\sup_{P\in \mathscr{P}}E_P[\varphi(X_1,\ldots,X_d)], \varphi\in C_b(\mathbb R^d). $$
Now, for each $l$, $\sum_{k=1}^l f_k(\bm x)$ is a bounded continuous function on $\mathbb R^{n_{l+1}}$. Thus
$$ \inf_{P\in \mathscr{P}}E_P[\sum_{k=1}^lf_k(\bm X)]=\inf_{P\in \mathrm{core}(\Sbep)}E_P[\sum_{k=1}^lf_k(\bm X)]\le E_Q[\sum_{k=1}^lf_k(\bm X)]\le M. $$
It follows that there exists $P_l\in \mathscr{P}$ such that
$$ E_{P_l}[\sum_{k=1}^lf_k(\bm X)]\le M+1/l. $$
By the weak compactness of $\mathscr{P}$, there exists a subsequence $l^{\prime}\nearrow \infty$ and $P\in \mathscr{P}$ such that $P_{l^{\prime}}$ converges weakly to $P$. It follows that
$$ E_{P}[\sum_{k=1}^lf_k(\bm X)]=\lim_{l^{\prime}\to \infty} E_{P_{l^{\prime}}}[\sum_{k=1}^lf_k(\bm X)]
\le \liminf_{l^{\prime}\to \infty}E_{P_{l^{\prime}}}[\sum_{k=1}^{l^{\prime}}f_k(\bm X)]\le M. $$
Hence
\begin{align} \label{eqconvergence3}
&\sum_{k=1}^{\infty} P\left(\max_{ n\le n_{k+1}}\left|\sum_{i=1}^n(X_i-\varphi_{i-1}(X_1,\ldots,X_{i-1}))\right|\ge 2\epsilon_k n_k\right)\nonumber \\
& \;\; \le E_{P}[\sum_{k=1}^{\infty}f_k(\bm X)]\le M<\infty,
\end{align}
where $M$ and $\{\epsilon_k\}$ do not depend on the $\varphi_i$'s and the choice of $P$. Let $\delta_n=\epsilon_k^{1/2}$ for $n_k+1\le n\le n_{k+1}$. Then \eqref{eqconvergence3} implies that
$$ P\left(\frac{\sum_{i=1}^n (X_i-\varphi_{i-1}(X_1,\ldots,X_{i-1}))}{\delta_n n}\to 0\right)=1, $$
and then
\begin{equation} \label{eqconvergence4}
P\left(\frac{\sum_{i=1}^n (X_i-\varphi_i(X_1,\ldots,X_i))}{\delta_n n}\to 0\right)=1.
\end{equation}
The proof is completed.

{\em Step 4.} We show (c). Let $\{\tilde{X}_n;n\ge 1\}$ be a copy on $\widetilde{\Omega}=\mathbb R^{\mathbb N}$ of $\{X_n;n\ge 1\}$, with $(\widetilde{\Omega},C(\widetilde{\Omega}),\widetilde{\mathbb E})$ and $\widetilde{\mathscr{P}}$ being defined as in Lemma \ref{lem:2}. Since condition (CC) is satisfied for $\widetilde{\mathbb E}$ with $\widetilde{\mathscr{P}}$,
by \eqref{eqconvergence4} there exists $\delta_k\searrow 0$ such that for every sequence of continuous functions $\varphi_i(x_1,\ldots,x_i)\in[\underline{\mu},\overline{\mu}]$ we have a $P\in \widetilde{\mathscr{P}}$ satisfying
$$ P\left(\frac{\widetilde{T}_n}{{\delta_n n}}\to 0\right)=1, $$
where $\widetilde{T}_n=\sum_{i=1}^n (\widetilde{X}_i-\varphi_i(\widetilde{X}_1,\ldots,\widetilde{X}_i))$,
which implies that
$$ \lim_{n\to \infty} \min_{N\ge n} P\left(\bigcap_{m=n}^N \left\{\frac{|\widetilde{T}_m|}{\delta_m m}\le 1/2\right\}\right)
=P\left( \frac{|\widetilde{T}_n|}{\delta_n n}\le 1/2\; i.o.\right)=1. $$
Let $ f(t)\in C_{b,Lip}(\mathbb R)$ be a function such that $I\{t\le 1\}\ge f(t)\ge I\{t\le 1/2\}$. We have
\begin{align*}
\min_{N\ge n}\sup_{P\in \widetilde{\mathscr{P}}} E_P f\left(\max_{n\le m\le N} \frac{|\widetilde{T}_m|}{\delta_m m} \right)
\ge \min_{N\ge n} P\left(\bigcap_{m=n}^N \left\{\frac{|\widetilde{T}_m|}{\delta_m m}\le 1/2\right\}\right)
\to 1.
\end{align*}
Let
$$ g(\bm x)=g(x_1,\ldots, x_N)=\max_{n\le m\le N}\frac{|\sum_{i=1}^m (x_i-\varphi_i(x_1,\ldots,x_i))|}{\delta_m m} \;\text{ and } $$
$$\varphi(\bm x)=\varphi(x_1,\ldots, x_N)=
f\left(\max_{n\le m\le N}\frac{|\sum_{i=1}^m (x_i-\varphi_i(x_1,\ldots,x_i))|}{\delta_m m}\right). $$
For \eqref{eqth2.6}, it is sufficient to show that
\begin{equation}\label{eq:proof:th2-c}
\upCapc\left(g(\bm X)\le 1\right) \ge \mathbb E^{\widetilde{\mathscr{P}}} \left[\varphi(\widetilde{\bm X})\right].
\end{equation}
Notice that $g$ and $\varphi$ are continuous functions on $\mathbb R^N$. But $\Sbep[\varphi(\bm X)]$ may have no definition. We need to approximate $\varphi$ by a bounded Lipschitz function.
Similar to \eqref{eq:finiteapprox}, for any $\epsilon>0$, there are a compact set $K=[-c,c]^{\bigotimes N}$ and $\psi(\bm x)=\psi(x_1,\ldots,x_N)
\in C_b^{\infty}(\mathbb R^N)$ such that
$$\sup_{\bm x}|\varphi(\bm x)-\psi(\bm x)|\le \frac{\epsilon}{2}+2I\{\bm x\in K^c\},
$$
$$\upCapc(\bm X\in K^c)\le \sum_{i=1}^N \upCapc(|X_i|>c)<\epsilon/4 \quad \text{and}\quad
\widetilde{\Capc}(\widetilde{\bm X}\in K^c)\le \sum_{i=1}^N \widetilde{\Capc}(|\widetilde{X}_i|>c)<\epsilon/4. $$
Since $\psi(x_1,\ldots,x_N)\in C_{b,Lip}(\mathbb R^N)$ and $(X_1,\ldots, X_N)\overset{d}
=(\widetilde{X}_1,\ldots, \widetilde{X}_N)$, we have
\begin{align*}
&\Sbep[\psi(\bm X)]=\widetilde{\mathbb E}[\psi(\widetilde{\bm X})]=\mathbb E^{\widetilde{\mathscr{P}}}[\psi(\widetilde{\bm X})]\\
\ge & \mathbb E^{\widetilde{\mathscr{P}}}[\varphi(\widetilde{\bm X})]-\epsilon/2-2\widetilde{\mathbb V}(\widetilde{\bm X}\in K^c)\ge
\mathbb E^{\widetilde{\mathscr{P}}}[\varphi(\widetilde{\bm X})]-\epsilon.
\end{align*}
On the other hand,
$$I\left\{g(\bm X)\le 1\right\}\ge \varphi(\bm X)\ge \psi(\bm X)-\epsilon/2-2I\{\bm X\in K^c\}. $$
It follows that
\begin{align*}
\upCapc\left(g(\bm X)\le 1\right)
\ge \Sbep[\psi(\bm X)]- \epsilon/2-2\upCapc(\bm X\in K^c) \ge
\mathbb E^{\widetilde{\mathscr{P}}}[\varphi(\widetilde{\bm X})]-2\epsilon.
\end{align*}
Thus, \eqref{eq:proof:th2-c} is proved.

Finally, by (3.5) in Lemma 3.1 of Zhang \cite{ZhangLIL21}, for any $\epsilon\in (0,1)$,
\begin{align*}
&\upCapc\left( \max_{n\le n_{k+1}}(S_n-n\overline{\mu})\ge \{\epsilon+2|\overline{\mu}-\Sbep[X_1^{(n_k)}]|\}n_k\right)\\
\le & \sum_{i=1}^{n_{k+1}}\upCapc\left(|X_i|>n_k\right)+ \upCapc\left( \max_{n\le n_{k+1}}\sum_{i=1}^n(X_i^{(n_k)}-\Sbep[X_i^{(n_k)}])\ge \epsilon n_k\right)\\
\le & \sum_{i=1}^{n_{k+1}}\upCapc\left(|X_i|\wedge n_k>n_k\right)+C
\frac{\sum_{i=1}^{n_{k+1}}\Sbep[(X_i^{(n_k)})^2]}{\epsilon^2n_k^2} \\
\le & \frac{C}{\epsilon^2} \frac{\Sbep[X_1^2\wedge n_k^2]}{n_k}
\le \frac{C}{\epsilon^2} \frac{\int_0^{n_k^2}\upCapc(X_1^2> x)dx}{n_k}
=\frac{2C}{\epsilon^2} \frac{\int_0^{n_k}y\upCapc(|X_1|>y)dy}{n_k}.
\end{align*}
Notice that
\begin{align*}
&\sum_{k=1}^{\infty} \frac{\int_0^{n_k}y\upCapc(|X_1|>y)dy}{n_k}\\
=&\int_0^1y\upCapc(|X_1|>y)dy\sum_{k=1}^{\infty} \frac{1}{n_k}
+\int_1^{\infty} y\upCapc(|X_1|>y)dy \sum_{k=1}^{\infty} \frac{I\{n_k\ge y\}}{n_k} \\
\le & \int_0^1\upCapc(|X_1|>y)dy +2 \int_1^{\infty} \upCapc(|X_1|>y)dy
\le 2C_{\upCapc}(|X_1|)< \infty,
\end{align*}
and
$$ |\overline{\mu}-\Sbep[X_1^{(n_k)}]|\le \vSbep[(|X_1|-n_k)^+]\to 0. $$
We have
$$ \sum_{k=1}^{\infty} \upCapc\left( \max_{n\le n_{k+1}}(S_n-n\overline{\mu})\ge \epsilon n_k\right)<\infty, \quad \forall \epsilon>0. $$
Similarly,
$$ \sum_{k=1}^{\infty} \upCapc\left( \max_{n\le n_{k+1}}(-S_n+n\underline{\mu})\ge \epsilon n_k\right)<\infty, \quad \forall \epsilon>0. $$
It follows that
$$\sum_{k=1}^{\infty} \upCapc\left( \max_{n\le n_{k+1}}(S_n-n\overline{\mu})\ge \epsilon n_k \text{ or } \max_{n\le n_{k+1}}(-S_n+n\underline{\mu})\ge \epsilon n_k\right)<\infty, \quad \forall \epsilon>0. $$
Thus, there exists a sequence $\epsilon_k\searrow 0$ such that
$$
\sum_{k=1}^{\infty} \upCapc \Big(\max_{ n\le n_{k+1}}\big(S_n- n\overline{\mu}\big)\ge \epsilon_k n_k \text{ or } \max_{ n\le n_{k+1}}\big(n\underline{\mu}-S_n\big)\ge \epsilon_k n_k\Big)<\infty,
$$
which implies \eqref{eqth2.4} and \eqref{eqth2.5}.
The proof of (c) is completed.
$\Box$

\bigskip
{\bf Proof of Theorem \ref{th3} and Remark \ref{rk:3}.} Let $\{\epsilon_i;i\ge 1\}$ satisfy \eqref{eqlem7.1}. For (a), let
$$\varphi_i(x_1,\ldots,x_i)=
\begin{cases}
0, & \text{if } i<d, \\
\big(\overline{\varphi}(x_1,\ldots,x_d)-\underline{\varphi}(x_1,\ldots,x_d)\big)\epsilon_i+ \underline{\varphi}(x_1,\ldots,x_d), &\text{if } i\ge d.
\end{cases} $$
Then by \eqref{eqlem7.1},
\begin{align*}
&\liminf\limits_{n\to \infty} \frac{\sum_{i=1}^n\varphi_i(x_1,\ldots,x_i)}{n}= \underline{\varphi}(x_1,\ldots,x_d),\\
&\limsup\limits_{n\to \infty} \frac{\sum_{i=1}^n\varphi_i(x_1,\ldots,x_i)}{n}= \overline{\varphi}(x_1,\ldots,x_d).
\end{align*}
Hence, by Theorem \ref{th2} (a) there is a probability measure $Q\in \mathrm{core}(\Sbep)$ such that
$$ Q\left(\liminf_{n\to \infty}\frac{S_n}{n} = \underline{\varphi}(X_1,\ldots,X_d) \text{ and } \limsup_{n\to \infty}\frac{S_n}{n} = \overline{\varphi}(X_1,\ldots,X_d)\right)=1,$$
which implies \eqref{eqth3.1}.

Next, we consider (b). For any $\epsilon>0$ and any $\varphi\in \mathbb L_c$ with $\varphi(\bm x)\in [\underline{\mu},\overline{\mu}]$, there is a bounded continuous function $f(\bm x)$ on $\mathbb R^{\mathbb N}$ such that $\mathbb E^{\mathscr{P}}[|\varphi(\bm X)-f(\bm X)|]< \epsilon/4$. Let $M=\sup_{\bm x}|f(\bm x)|$. Choose $c_k$ such that $\upCapc(|X_k|\ge c_k)< \epsilon /(2^k(4M))$. Let $K=\bigotimes_{k=1}^{\infty}[-c_k,c_k]$. Then $K$ is a compact set in $\mathbb R^{\mathbb N}$. For $\bm x=(x_1,x_2,\ldots)$, denote $\bm x^d=(x_1,\ldots,x_d,x_{d+1}^0, x_{d+2}^0,\ldots)$ and $f_d(\bm x)=f(\bm x^d)$. Then $f_d(\bm x)$ is a bounded continuous function on $\mathbb R^d$, and
$$ \sup_{\bm x}|f_d(\bm x)-f(\bm x)|\le \sup_{\bm x\in K}|f(\bm x^d)-f(\bm x)|+2MI\{\bm x\in K^c\}. $$
Notice that $f(\bm x)$ is uniformly continuous on the compact set $K$, and $d(\bm x^d,\bm x)\le \sum_{j=d+1}^{\infty}1/2^j=1/2^d\to 0$. We have, for $d$ large enough,
$$\sup_{\bm x\in K}|f(\bm x^d)-f(\bm x)| <\epsilon/4.$$
It follows that
$$\mathbb E^{\mathscr{P}}[|f_d(\bm X)-f(\bm X)|]\le \epsilon/4+2M\sum_{k=1}^{\infty}\upCapc(|X_k|\ge c_k)\le \epsilon/4+\epsilon/2. $$
Therefore, for any $\epsilon>0$, there exist $d$ and a bounded continuous function $f_d(\bm x)$ on $\mathbb R^d$ such that
$$ \mathbb E^{\mathscr{P}}[|\varphi(\bm X)-f_d(\bm X)|]\le\Sbep[|\varphi(\bm X)-f(\bm X)|]+\Sbep[|f(\bm X)-f_d(\bm X)|]<\epsilon. $$
Hence, there exists a sequence $d_l\nearrow \infty$ and bounded continuous functions $f_{d_l}$ on $\mathbb R^{d_l}$ such that
$\sum_{l=1}^{\infty} \mathbb E^{\mathscr{P}}[|\varphi(\bm X)-f_{d_l}(\bm X)|]<\infty$. Let
$\Omega_0=\{\sum_{l=1}^{\infty}|\varphi(\bm X)-f_{d_l}(\bm X)|<\infty\}$.
Then
$$
f_{d_l}(\bm X)\to \varphi(\bm X) \text{ on } \Omega_0.
$$
For each $P\in\mathscr{P}$, $\sum_{k=1}^{\infty}E_P[|\varphi(\bm X)-f_{d_l}(\bm X)|]\le\sum_{l=1}^{\infty} \mathbb E^{\mathscr{P}}[|\varphi(\bm X)-f_{d_l}(\bm X)|]<\infty$, which implies $P(\Omega_0^c)=0$.
Hence $\mathbb V^{\mathscr{P}}(\Omega_0^c)=0$.
Let $d_0=0$ and $f_{d_0}$ be a constant in $[\underline{\mu},\overline{\mu}]$. Denote $l_i=\max\{l:d_{l}\le i\}$,
$$ \varphi_i(x_1,\ldots,x_i)=\underline{\mu}\vee f_{d_{l_i}}(\bm x)\wedge \overline{\mu}. $$
Then $\varphi_i(x_1,\ldots,x_i)$ is a continuous function on $\mathbb R^i$ and
$$ \mathbb V^{\mathscr{P}}\left(\lim_{n\to \infty} \frac{\sum_{i=1}^n \varphi_i(X_1,\ldots,X_i)}{n}\ne \varphi(\bm X) \right)=0. $$
By Theorem \ref{th2} (b), there exists $P_{\varphi}\in\mathscr{P}$ such that \eqref{eqcor1.3} holds with $\varphi$.

Now, for $\underline{\varphi}, \overline{\varphi}\in \mathbb L_c$ with $\underline{\mu}\le \underline{\varphi}\le \overline{\varphi}\le \overline{\mu}$, there are continuous functions $\underline{\varphi}_i(x_1,\ldots,x_i)$, $\overline{\varphi}_i(x_1,\ldots,x_i)$, $i=1,2,\ldots$, and an event $\Omega_0$ such that $\mathbb V^{\mathscr{P}}(\Omega_0^c)=0$,
\begin{equation}\label{eqproofth3.3}
\lim_{i\to \infty}\underline{\varphi}_i(X_1,\ldots,X_i) = \underline{\varphi}(\bm X) \quad \text{and}\quad \lim_{i\to \infty}\overline{\varphi}_i(X_1,\ldots,X_i) = \overline{\varphi}(\bm X)
\quad \text{on } \Omega_0.
\end{equation}
Without loss of generality, we can assume $\underline{\mu}\le \underline{\varphi}_i(x_1,\ldots,x_i)\le \overline{\varphi}_i(x_1,\ldots,x_i)\le \overline{\mu}$. Let
$$\varphi_i(x_1,\ldots,x_i)=
\big(\overline{\varphi}_i(x_1,\ldots,x_i)-\underline{\varphi}_i(x_1,\ldots,x_i)\big)\epsilon_i+ \underline{\varphi}_i(x_1,\ldots,x_i). $$
Then on the event $\Omega_0$,
\begin{align*}
&\liminf\limits_{n\to \infty} \frac{\sum_{i=1}^n\varphi_i(X_1,\ldots,X_i)}{n}= \underline{\varphi}(\bm X),\\
&\limsup\limits_{n\to \infty} \frac{\sum_{i=1}^n\varphi_i(X_1,\ldots,X_i)}{n}= \overline{\varphi}(\bm X),
\end{align*}
by \eqref{eqlem7.1} and \eqref{eqproofth3.3}.
Hence, by Theorem \ref{th2} (b) there is a probability measure $P\in\mathscr{P}$ such that
$$ P\left(\liminf_{n\to \infty}\frac{S_n}{n} = \underline{\varphi}(\bm X) \text{ and } \limsup_{n\to \infty}\frac{S_n}{n} = \overline{\varphi}(\bm X)\right)=1,$$
which implies \eqref{eqth3.2}. $\Box$

\bigskip
{\bf Proof of Corollary \ref{cor2}.} \eqref{eqcor2.1} and \eqref{eqcor2.2} follow immediately from \eqref{eqcor1.2} and \eqref{eqcor1.3}, respectively.

For (c), \eqref{eqcor2.3} was proved by Zhang \cite{ZhangLLN}. For \eqref{eqcor2.4}, we consider the copy $\{\tilde{X}_n; n\ge 1\}$ on $\mathbb R^{\mathbb N}$ of $\{X_n;n\ge 1\}$ as defined in Lemma \ref{lem:2}. By Corollary \ref{cor2} (b), for any continuous function $\varphi(\bm x):\mathbb R^{\mathbb N}\to [\underline{\mu},\overline{\mu}]$, there is a probability measure $\widetilde{P}_{\varphi}\in \widetilde{\mathscr{P}}_{\max}$ such that
\begin{align}\label{eq:proofcor2.7}
&\widetilde{\cCapc}\left(\Big|\frac{\sum_{i=1}^n \tilde{X}_i}{n}-\varphi(\tilde{\bm X})\Big|\ge \epsilon \right)\nonumber \\
\le & \widetilde{P}_{\varphi}\left(\Big|\frac{\sum_{i=1}^n \tilde{X}_i}{n}-\varphi(\tilde{\bm X})\Big|\ge \epsilon \right)\to 0 \text{ for all } \epsilon>0.
\end{align}
Notice that $\varphi(\bm X)$ may not be an element of $\mathscr{H}$. We cannot use inequality \eqref{eqV-V} directly. We will consider a finite-dimensional approximation of $\varphi(\bm X)$.
Let $M=\sup_{\bm x}|\varphi(\bm x)|$, and let $K$ be a compact set. As shown in the proof of Lemma \ref{lem:1} (cf. \eqref{eq:finiteapprox}), there exists $\psi\in C_{b}^{\infty}(\mathbb R^d)$ such that
\begin{equation}\label{eq:proofcor2.8}
\sup_{\bm x}|\varphi(\bm x)-\psi(\bm x)|\le \frac{\epsilon}{2}I\{\bm x\in K\}+2MI\{\bm x\in K^c\}.
\end{equation}
It follows that
$$\cCapc\left(\Big|\frac{S_n}{n}-\varphi(\bm X)\Big|\ge 5\epsilon \right)
\le \cCapc\left(\Big|\frac{S_n}{n}-\psi(\bm X)\Big|\ge 4\epsilon \right)+\Capc(\bm X\in K^c),
$$
$$\widetilde{P}_{\varphi}\left(\Big|\frac{\sum_{i=1}^n \tilde{X}_i}{n}-\psi(\tilde{\bm X})\Big|\ge 2\epsilon \right)
\le \widetilde{P}_{\varphi}\left(\Big|\frac{\sum_{i=1}^n \tilde{X}_i}{n}-\varphi(\tilde{\bm X})\Big|\ge \epsilon \right)
+\widetilde{P}_{\varphi}(K^c). $$
Let $f_{\epsilon}\in C_{b,Lip}(\mathbb R)$ such that $I\{|x|\ge 2\epsilon\}\le f_{\epsilon}(x)\le I\{|x|\ge \epsilon\}$. Notice
$f_{\epsilon}\Big(\frac{\sum_{i=1}^n x_i}{n}-\psi(\bm x)\Big)\in C_{b,Lip}(\mathbb R^{n+d})$ and $f_{\epsilon}\Big(\frac{S_n}{n}-\psi(\bm X)\Big)\in \mathscr{H}$. Then
\begin{align*}
\cCapc&\left(\Big|\frac{S_n}{n}-\psi(\bm X)\Big|\ge 4\epsilon \right)
\le \cSbep\left[f_{2\epsilon}\Big(\frac{S_n}{n}-\psi(\bm X)\Big)\right] \; (\text{by } \eqref{eq1.4})  \\
=& \inf_{P\in \widetilde{\mathscr{P}}_{\max}}E_P\left[f_{2\epsilon}\Big(\frac{\sum_{i=1}^n \tilde{X}_i}{n}-\psi(\tilde{\bm X})\Big)\right] \; (\text{by } \eqref{eq:lem2.0}) \\
\le &E_{\widetilde{P}_{\varphi}} \left[f_{2\epsilon}\Big(\frac{\sum_{i=1}^n \tilde{X}_i}{n}-\psi(\tilde{\bm X})\Big)\right]\le \widetilde{P}_{\varphi}\left(\Big|\frac{\sum_{i=1}^n \tilde{X}_i}{n}-\psi(\tilde{\bm X})\Big|\ge 2\epsilon \right).
\end{align*}
Hence
\begin{align*}
& \cCapc\left(\Big|\frac{S_n}{n}-\varphi(\bm X)\Big|\ge 5\epsilon \right) \\
\le & \widetilde{P}_{\varphi}\left(\Big|\frac{\sum_{i=1}^n \tilde{X}_i}{n}-\varphi(\tilde{\bm X})\Big|\ge \epsilon \right)
+\Capc(\bm X\in K^c)+ \widetilde{P}_{\varphi}(K^c).
\end{align*}
If $\varphi(x_1,\ldots,x_d)$ is a continuous function on $\mathbb R^d$, then $K=\bigoplus_{k=1}^d[-l2^k,l2^k]$ is a compact set in $\mathbb R^d$. By the (finite) sub-additivity of $\Capc$, we have
$$ \Capc(\bm X\in K^c)\le \sum_{k=1}^d\Capc(|X_k|>l2^k)\le \sum_{k=1}^d\frac{\vSbep[|X_1|]}{l2^k}\le \frac{\vSbep[|X_1|]}{l}
\to 0  $$
as $l\to \infty$. Similarly, $\widetilde{P}_{\varphi}(K^c)\le \frac{\vSbep[|X_1|]}{l}\to 0$ as $l\to \infty$. Thus,
\begin{equation}\label{eq:proofcor2.9}
\cCapc\left(\Big|\frac{S_n}{n}-\varphi(\bm X)\Big|\ge 5\epsilon \right)
\le \widetilde{P}_{\varphi}\left(\Big|\frac{\sum_{i=1}^n \tilde{X}_i}{n}-\varphi(\tilde{\bm X})\Big|\ge \epsilon \right)\to 0,
\end{equation}
as $n\to \infty$, by \eqref{eq:proofcor2.7}.

If $\varphi(\bm x)$ is a uniformly continuous function on $\mathbb R^{\mathbb N}$, then \eqref{eq:proofcor2.8} holds with $K=\mathbb R^{\mathbb N}$. And so, \eqref{eq:proofcor2.9} remains true.

If $\varphi(\bm x)$ is a continuous function on $\mathbb R^{\mathbb N}$ and the capacity is countably sub-additive, then we choose the compact set $K=\bigoplus_{k=1}^{\infty}[-l2^k,l2^k]$ in $\mathbb R^{\mathbb N}$. By the countable sub-additivity of $\Capc$, we have
$$ \Capc(\bm X\in K^c)\le \sum_{k=1}^{\infty}\Capc(|X_k|>l2^k)\le \sum_{k=1}^{\infty}\frac{\vSbep[|X_1|]}{l2^k}\le \frac{\vSbep[|X_1|]}{l}
\to 0,  $$
and similarly, $\widetilde{P}_{\varphi}(K^c)\le \frac{\vSbep[|X_1|]}{l}\to 0$ as $l\to \infty$.
\eqref{eq:proofcor2.9} remains true. $\Box$.

\bigskip
{\bf Proof of Corollary \ref{rk:5}.} We can consider $\mathscr{H}(\bm X)$ instead of $\mathscr{H}$. By Lemma \ref{lem:0} (a), it is sufficient to show that $\Sbep$ is regular on $\mathscr{H}(\bm X)$.
Let $\mathscr{H}(\bm X)\ni Z_n\searrow 0$. We need to show that $\Sbep[Z_n]\to 0$. For each $Z_n$, by definition there is a $\varphi_n(x_1,\ldots,x_{d_n})\in C_{b,Lip}(\mathbb R^{d_n})$ such that $Z_n=\varphi_n(X_1,\ldots, X_{d_n})$. For any $\epsilon>0$, let $\Omega_n=\{Z_n\ge \epsilon\}$, $F_n=\{\bm x; \varphi_n(x_1,\ldots,x_{d_n})\ge \epsilon\}$, $C_n=\bigcap_{m=1}^n F_m$, and $C=\bigcap_{m=1}^{\infty} F_m$. Then $F_n,C_n,C$ are closed sets of $\mathbb R^{\mathbb N}$, $\Omega_n\searrow\emptyset$, $C_n\searrow C$,
$$ \Omega_n=\bigcap_{m=1}^n \Omega_m= \bigcap_{m=1}^n\{\bm X\in F_m\} =\{\bm X\in C_n\}  $$
and
$$ \{\bm X\in C\}= \bigcap_{m=1}^{\infty}\{\bm X\in F_m\} =\lim_n\{Z_n\ge \epsilon\}=\emptyset.$$
By the assumption that $\Capc\bm X^{-1}$ is continuous from above on closed sets, we have
\begin{align*}
0\le \Sbep[Z_n]\le& \epsilon+\Capc(Z_n\ge \epsilon)=\epsilon+\Capc(\bm X\in C_n)\\
&\searrow \epsilon+\Capc(\bm X\in C)=\epsilon.
\end{align*}
Thus, $\Sbep[Z_n]\to 0$. The proof is completed. $\Box$.

\bigskip
{\bf Proof of Proposition \ref{th:example}.} Consider the subset $\mathbb X=\{0,1,2\}^{\mathbb N}$ of $\mathbb R^{\mathbb N}$. Then $\mathbb X$ is a compact set of $\mathbb R^{\mathbb N}$. We choose $\mathscr{P}$ and $f$ as follows:
$$ \mathscr{P}=\{\delta_{\bm x}:\bm x\in \mathbb X\}, \quad A=\{\bm x=(x_1,x_2,\ldots)\in \mathbb X: x_n=2 \text{ i.o.}\}, \quad f=I_A, $$
where $\delta_{\bm x}$ is the unit mass at $\bm x$, i.e., $\delta_{\bm x}(B)=I_B(\bm x)$.

First, we show that $\mathscr{P}$ is a weakly compact family of probabilities on $\mathbb R^{\mathbb N}$. Suppose $P_n=\delta_{\bm x_n}\in \mathscr{P}$. Notice that $\mathbb X$ is compact. There are a subsequence $n_k$ and a point $\bm x\in \mathbb X$ such that $\bm x_{n_k}\to \bm x$. Thus, for any continuous function $h$ on $\mathbb R^{\mathbb N}$,
$$ \int h\, d\, P_{n_k}= h(\bm x_{n_k}) \to h(\bm x)=\int h \,d \,\delta_{\bm x}. $$
Hence $P_{n_k}\implies \delta_{\bm x}\in \mathscr{P}$. It follows that $\mathscr{P}$ is a weakly compact family.

Second, we show that $f$ is a Borel function. It is obvious that $\{\bm x: x_n=2\}$ is a closed set in $\mathbb R^{\mathbb N}$. Hence
$$ A=\mathbb X\bigcap \bigcap_{m=1}^{\infty}\bigcup_{n=m}^{\infty} \{\bm x: x_n=2\}$$
is a Borel set, and thus $f=I_A$ is a bounded Borel function.

Next, suppose $\mathbb E^{\mathscr{P}} |f_n-f|\to 0$, where the $f_n$'s are finite-dimensional Borel functions. Without loss of generality, it can be assumed that $f_n(\bm x)=f_n(x_1,\ldots,x_n)$. For any $\bm x\in \mathbb X$,
$$ |f_n(\bm x)-f(\bm x)|=\int |f_n-f|d\, \delta_{\bm x}\le \mathbb E^{\mathscr{P}} |f_n-f|\to 0. $$
That is
$$ f_n\big|_{\mathbb X}\to f\big|_{\mathbb X}. $$
On the other hand, $f_n$ is a continuous function on $\{0,1,2\}^n$ since $\{0,1,2\}^n$ has only finitely many points, and so a continuous function on $\mathbb X$. It follows that, on the complete metric space $\mathbb X$, $f$ is a pointwise limit of a sequence of continuous functions, and thus its continuous points are dense in $\mathbb X$ by the Baire--Osgood theorem. On the other hand, it is easily seen that both $A$ and $\mathbb X\setminus A$ are dense in $\mathbb X$ since $A\ni(x_1,\ldots,x_n,2,2,\ldots)\to \bm x$ and $\mathbb X\setminus A \ni (x_1,\ldots,x_n,0,0,\ldots)\to \bm x$ as $n\to\infty$ for any $\bm x=(x_1,x_2,\ldots)\in \mathbb X$, which implies that $f=I_A$ is nowhere continuous on $\mathbb X$. We arrive at a contradiction.

Finally, notice
$$ \mathbb E^{\mathscr{P}}\left[\varphi(X_1,\ldots, X_d)\right]=
\sup_{(x_1,\ldots, x_d)\in\{0,1,2\}^d} \varphi(x_1,\ldots, x_d). $$
It is easily seen that $\{X_n;n\ge 1\}$ is a sequence of i.i.d. random variables under $\mathbb E^{\mathscr{P}}$. See the proof of Theorem \ref{th4} below.
$\Box$

\bigskip
{\bf Proof of Theorem \ref{th4}.} Instead of $\Omega=\mathbb N\cup\{0\}$ as in Ter\'an \cite{Teran18}, we consider the space $\Omega=[0,\infty)$, take $\mathscr{H}$ to be the set of all bounded Borel-measurable functions on $\Omega$, and define
$\Sbep[X] = \sup_{\omega\in\Omega} X(\omega)$ for all $X\in \mathscr{H}$. It is easily checked that $\Sbep$ is a sub-linear expectation. For this sub-linear expectation space, $\Sbep[X]=\sup_{P\in \mathscr{P}}E_P[X]$, where $\mathscr{P}=\{\delta_{\omega_0};\omega_0\in \Omega\}$, $\delta_{\omega_0}$ is the unit mass at $\omega_0$, i.e., $\delta_{\omega_0}(A)=I_A(\omega_0)$. Also, $\mathrm{core}(\Sbep)=\{\text{all probability measures on } \Omega\}$, $\upCapc(A)=\outCapc(A)=\Capc^{\mathscr{P}}(A)=\Sbep[I_A]=\sup_{\omega\in \Omega}I_A(\omega)=1$ if $A$ is not empty and $0$ otherwise,
$\lowCapc(A)=\outcCapc(A)=\cCapc^{\mathscr{P}}(A)=\cSbep[I_A]=\inf_{\omega\in \Omega}I_A(\omega)=1$ if $A=\Omega$ and $0$ otherwise. This sub-linear expectation $\Sbep$ is not regular because $Y_n(\omega)=\frac{\omega\wedge n}{n}\searrow 0$, but
$\Sbep[Y_n]=\sup_{\omega}Y_n(\omega)=1\not\to 0$.

Let $\epsilon_n(k)$ be the $n$th bit in the ternary representation of a positive integer $k\in \mathbb N\cup\{0\}$, and $\zeta_n(\omega_0)$ be the $n$th bit in the ternary representation of a decimal $\omega_0\in [0,1)$, respectively, i.e.
$$ k=\sum_{n=1}^{\infty}3^{n-1}\epsilon_n(k), \; k\in \mathbb N\cup\{0\}, \quad \omega_0 =\sum_{n=1}^{\infty} 3^{-n}\zeta_n(\omega_0), \; \omega_0\in [0,1), $$
with $\epsilon_n(k), \zeta_n(\omega_0)\in \{0, 1,2\}$.
Write $\bm \epsilon(k)=(\epsilon_1(k),\epsilon_2(k),\ldots)$ and $\bm \zeta(\omega_0)=(\zeta_1(\omega_0),\zeta_2(\omega_0),\ldots)$.
For $\omega\in\Omega$, we define $\bm\epsilon(\omega)=\bm\epsilon([\omega])$ and $\bm\zeta(\omega)=\bm\zeta(\omega-[\omega])$,
where $[\omega]$ is the integer part of $\omega$. Then $\omega$ has a unique representation:
\begin{equation}\label{eq:ternary}
\omega=\sum_{n=1}^{\infty}3^{n-1}\epsilon_n(\omega)+\sum_{n=1}^{\infty} 3^{-n}\zeta_n(\omega), \quad \omega\in \Omega=[0,\infty).
\end{equation}
Notice
$$ f(\bm x) =\sum_{k=1}^{n}3^{k-1}x_k: \{0,1,2\}^n  \to \{0,1,\ldots, 3^n-1\} $$
$$ g(\bm x)=\sum_{k=1}^{\infty}3^{-k}x_k: \{0,1,2\}^{\mathbb N}\setminus \{  (2,2,\ldots)\}\to [0,1) $$
are both one-to-one maps.
Hence
\begin{align*}
&  \{(\epsilon_1(k),\ldots,\epsilon_n(k)):k\in\mathbb N\cup\{0\}\}=\{0, 1,2\}^n, \\
&   \{(\zeta_1(\omega_0),\ldots,\zeta_n(\omega_0)):\omega_0\in [0,1)\}=\{0, 1,2\}^n.
\end{align*}
Thus
\begin{align*}
&\{(\epsilon_1(\omega),\ldots,\epsilon_n(\omega), \zeta_1(\omega),\ldots,\zeta_m(\omega)):\omega\in\Omega\}\\
=& \{(\epsilon_1(k),\ldots,\epsilon_n(k), \zeta_1(\omega_0),\ldots,\zeta_m(\omega_0)):k\in\mathbb N, \omega_0\in (0,1)\}
=\{0, 1,2\}^{n+m},
\end{align*}
and
$$ \{(\eta_1(\omega),\ldots,\eta_{n+m}(\omega)):\omega\in\Omega\}=\{0, 1,2\}^{n+m}, $$
if $\{\eta_1,\ldots,\eta_{n+m}\}$ is a permutation of $\{\epsilon_1,\ldots,\epsilon_n,\zeta_1,\ldots,\zeta_m\}$.
With the same arguments as Ter\'an \cite{Teran18}, it can be shown that $\{\eta_1,\ldots,\eta_{n+m}\}$ are i.i.d. random variables under $\Sbep$. In fact, for any $\varphi(x)\in C_{b,Lip}(\mathbb R^d)$ $(d\le n+m)$,
$$ \varphi_{d-1}(x_1,\ldots,x_{d-1})=:\Sbep[\varphi(x_1,\ldots,x_{d-1},\eta_d)]
=\sup_{x_d=0,1,2}\varphi(x_1,\ldots,x_{d-1},x_d), $$
and thus
\begin{align*}
&\Sbep\left[\Sbep[\varphi(x_1,\ldots,x_{d-1},\eta_d)]\big|_{x_1=\eta_1,\ldots,\eta_{d-1}}\right]
= \Sbep[\varphi_{d-1}(\eta_1,\ldots,\eta_{d-1})] \\
=&\sup_{\omega\in\Omega}\varphi_{d-1}\big(\eta_1(\omega),\ldots,\eta_{d-1}(\omega)\big)
=\sup_{(x_1,\ldots,x_{d-1})\in \{0,1,2\}^{d-1}}\varphi_{d-1}(x_1,\ldots,x_{d-1})\\
=&\sup_{(x_1,\ldots,x_{d-1})\in \{0,1,2\}^{d-1}}[\sup_{x_d=0,1,2}\varphi(x_1,\ldots,x_{d-1},x_d)]\\
=&\sup_{(x_1,\ldots,x_d)\in \{0,1,2\}^d} \varphi(x_1,\ldots, x_d) =\sup_{\omega\in\Omega}\varphi\big(\eta_1(\omega),\ldots,\eta_{d}(\omega)\big)
= \Sbep[\varphi(\eta_1,\ldots,\eta_d)],
\end{align*}
which implies independence. The identical distribution is obvious.

Let $h(x):\{0, 1,2\}\to \mathbb R$ be a function such that
$$ h(x)= \mu \text{ if } x= 0, \quad
\underline{\mu} \text{ if } x= 1,\quad
\overline{\mu} \text{ if } x=2,
$$
and $X_n(\omega)=h(\epsilon_n(\omega))$, $Y_n(\omega)=h(\zeta_n(\omega))$. Then $\{X_n;n\ge 1\}$ and $\{Y_n;n\ge 1\}$ satisfy the properties in (a) with
$$\Sbep[\varphi(X_n)]=\Sbep[\varphi(Y_n)]
=\max_{x=0,1,2}\varphi(h(x))=\max\{\varphi(\underline{\mu}), \varphi(\mu),\varphi(\overline{\mu})\},$$
$$\cSbep[\varphi(X_n)]=\cSbep[\varphi(Y_n)]
=\min_{x=0,1,2}\varphi(h(x))=\min\{\varphi(\underline{\mu}), \varphi(\mu),\varphi(\overline{\mu})\}.$$
Theorem \ref{th4} (a) and (b) are verified.

By \eqref{eq:ternary}, for each $\omega\in\Omega$, $\epsilon_i(\omega)=0$ when $3^{i-1}>\omega$. Hence
\begin{align*}
&\frac{1}{n}\sum_{i=1}^n g(X_i(\omega))=\frac{1}{n}\sum_{i=1}^n g\big(h(\epsilon_i(\omega))\big)\\
=&\frac{1}{n}\sum_{i: 3^{i-1}\le \omega}g\big(h(\epsilon_i(\omega))\big)
+\frac{1}{n}\sum_{i: \log_3\omega+1<i\le n}g(h(0))\\
&\to g(h(0))=g(\mu).
\end{align*}
\eqref{eq:th4.2} is proved. \eqref{eq:th4.0} and \eqref{eq:th4.1} are obvious. So, (c) is proved.

(d) follows from Corollary \ref{cor2} (c) by noting that $\upCapc(A)=\sup_{\omega\in \Omega}I_A(\omega)$ is a countably sub-additive capacity which satisfies \eqref{eq1.4}.

For (e), let $\{\tilde{X}_n;n\ge 1\}$ be a copy on $\{\underline{\mu},\mu,\overline{\mu}\}^{\mathbb N}$ of $\{X_n;n\ge 1\}$ (also a copy of $\{Y_n;n\ge 1\}$) defined as in Lemma \ref{lem:2}. By Theorem \ref{th3}, there is a $Q\in \widetilde{\mathscr{P}}_{\max}$ such that
$$ Q\left(C\Big\{\frac{\sum_{i=1}^n \tilde X_i}{n}\Big\}=\Big[\underline{\varphi}(\tilde{\bm X}),
\overline{\varphi}(\tilde{\bm X})\Big]\right)=1. $$
In the construction of $\widetilde{Q}$ (cf. Step 2 of the proof of Theorem \ref{th2}), we can choose $\widetilde{Q}_1$ such that $Q(\tilde{X}_1=\underline{\mu})=1$.
Hence
$$ Q\left(C\Big\{\frac{\sum_{i=1}^n \tilde X_i}{n}\Big\}=\Big[\underline{\varphi}(\tilde{\bm X}),
\overline{\varphi}(\tilde{\bm X})\Big], \tilde{X}_1=\underline{\mu}\right)=1. $$
Thus
$$ \left\{C\Big\{\frac{\sum_{i=1}^n \tilde x_i}{n}\Big\}=\Big[\underline{\varphi}(\tilde{\bm x}),
\overline{\varphi}(\tilde{\bm x})\Big], \tilde{x}_1=\underline{\mu}\right\} $$
is a non-empty set, and so
\begin{align*}\Delta_0=& \left\{\tilde{\bm \zeta}=(\tilde\zeta_1,\tilde\zeta_2,\ldots)\in \{0,1,2\}^{\mathbb N}
: \tilde\zeta_1=1, \right.\\
\quad  &\left.C\Big\{\frac{\sum_{i=1}^n h(\widetilde{\zeta}_i)}{n}\Big\}=\Big[\underline{\varphi}\big(h(\tilde{\zeta}_1), h(\tilde{\zeta}_2),\ldots\big),
\overline{\varphi}\big(h(\tilde{\zeta}_1), h(\tilde{\zeta}_2),\ldots\big)\Big]\right\}
\end{align*}
is a non-empty set. Since
$$ \omega_0=k+\sum_{n=1}^{\infty}3^{-n}\zeta_n(\omega_0): [k,k+1)\to \{0,1,2\}^{\mathbb N}\setminus\{ (2,2,\ldots)\} $$
is a one-to-one map, the set
$\left\{\omega\in \Omega: \bm\zeta(\omega)=(\zeta_1(\omega),\zeta_2(\omega),\ldots)\in \Delta_0\right\}$ is non-empty. Choose $\omega_0\in \left\{\omega\in \Omega: \bm\zeta(\omega)\in \Delta_0\right\}$, $P=\delta_{\omega_0}\in \mathscr{P}$.
It follows that
\begin{align*}
\upCapc \left(C\Big\{\frac{\sum_{i=1}^n Y_i}{n}\Big\}=\Big[\underline{\varphi}(\bm Y),\overline{\varphi}(\bm Y) \Big]\right)
\ge & P \left(C\Big\{\frac{\sum_{i=1}^n Y_i}{n}\Big\}=\Big[\underline{\varphi}(\bm Y),\overline{\varphi}(\bm Y) \Big]\right) \\
\ge &  P\left(\omega\in \Omega: \bm\zeta(\omega)\in \Delta_0\right)=1.
\end{align*}
The proof of (e) is completed.

For (f), let $\underline{\varphi}=\overline{\varphi}=\varphi$. By (e), there exists $P\in \mathscr{P}$ such that
$$ P\left(\frac{\sum_{i=1}^n Y_i}{n}\to \varphi(\bm Y)\right)=1, $$
which implies $\frac{\sum_{i=1}^n Y_i}{n}\overset{P}\to \varphi(\bm Y)$. The proof of (f) is completed.

Finally, we consider (g). Define the function
$$ h_n(\bm x)=h_n(x_1,\ldots,x_n)=\frac{1}{1\vee\big(n \sum_{k=1}^n \frac{ |x_k-\overline{\mu}|}{ 3^{k}}\big)}. $$
Then $h_n(\bm x)\in C_{b,Lip}(\mathbb R^{n})$, and $h_n(\bm x)\searrow 0$ if $\bm x\ne (\overline{\mu},\overline{\mu},\ldots)$ as $n\to\infty$, $h_n(\bm x)\equiv 1$ if $\bm x= (\overline{\mu},\overline{\mu},\ldots)$.

Since $\{(\epsilon_1(\omega),\ldots,\epsilon_n(\omega),\zeta_1(\omega),\ldots,\zeta_n(\omega)): \omega\in \Omega\}=\{0,1,2\}^{n+m}$ for all $n$, $m$, we have $\{(X_1(\omega),\ldots,X_n(\omega)): \omega\in \Omega\}=\{(Y_1(\omega),\ldots,Y_n(\omega)): \omega\in \Omega\}=\{\mu,\underline{\mu},\overline{\mu}\}^n$. Then
\begin{align*}
\Sbep[h_n(\bm X)] =\Sbep[h_n(\bm Y)]=\sup_{(x_1,\ldots, x_n)\in \{\mu,\underline{\mu},\overline{\mu}\}^n}h_n(\bm x) =h_n(\overline{\mu},\ldots,\overline{\mu})=1.
\end{align*}
On the other hand, it is obvious that $\bm\epsilon(\omega)\ne \bm (2,2,\ldots)$ and $\bm\zeta(\omega)\ne \bm (2,2,\ldots)$. Since $h(x): \{0,1,2\}\to \{\mu, \underline{\mu},\overline{\mu}\}$ is a one-to-one map, $\bm X(\omega)\ne (\overline{\mu},\overline{\mu},\ldots)$ and $\bm Y(\omega)\ne (\overline{\mu},\overline{\mu},\ldots)$.
So, $ h_n(\bm X(\omega))\searrow 0$ and $ h_n(\bm Y(\omega))\searrow 0$, $\omega\in\Omega$.
It follows that $\Sbep$ on $\mathscr{H}(\bm X)$ or $\mathscr{H}(\bm Y)$ is not regular.
The proof is completed. $\Box$

\begin{remark}
For the sub-linear expectation space $(\Omega,\mathscr{H},\Sbep)$ in Theorem \ref{th4}, $\cCapc=\lowCapc$ satisfies the following properties.
\begin{description}
\item[\rm (i)] $\cCapc(\emptyset)=0$, $\cCapc(\Omega)=1$.
\item[\rm (ii)] $A\subset B\implies \cCapc(A)\le \cCapc(B)$.
\item[\rm (iii)] $B_n\searrow B\implies \cCapc(B_n)\to \cCapc(B)$ whenever $B_n,B\in \mathcal{B}(\Omega)$.
\item[\rm (iv)] $\cCapc$ is completely monotone:
$$ \cCapc\left(\cup_{i=1}^n A_i\right)\ge \sum_{I\subset\{1,\ldots,n\}, I\ne \emptyset}(-1)^{|I|+1}\cCapc\left(\cap_{i\in I}A_i\right)\; \text{ for any } n\in \mathbb N. $$
\item[\rm (v)] $C_{\cCapc}(X)=\cSbep[X]=\inf_{\omega}X(\omega)$, $C_{\Capc}(X)=\Sbep[X]=\sup_{\omega}X(\omega)$.
\item[\rm (vi)] $\{\epsilon_n,\zeta_n;n\ge 1\}$ are independent (preindependent) and identically distributed relative to $\cCapc$ in the sense that
$$ \cCapc\left(\bigcap_{i=1}^k\{\epsilon_{n_i}\in A_i\},\bigcap_{j=1}^l\{\zeta_{m_j}\in B_j\}\right)
=\prod_{i=1}^k\cCapc\left(\epsilon_1\in A_i\right)\prod_{j=1}^l\cCapc\left(\epsilon_1\in B_j\right) $$
for any choice of $\{n_1, \ldots, n_k, m_1,\ldots,m_l\} \subset \mathbb N$, $A_i,B_j\in \mathcal{B}(\Omega)$, and
$k,l\in \mathbb N$.
\end{description}
Theorem \ref{th4} also gives different behaviours of the i.i.d. random variables under a capacity.
\end{remark}

\end{document}